\documentclass[preprint,12pt]{elsarticle}

\usepackage{amssymb}
\usepackage{amsmath}
\usepackage{mathtools}
\usepackage{xcolor}
\usepackage{subcaption}
\usepackage{multirow}
\usepackage{array}
\PassOptionsToPackage{hyphens}{url}
\usepackage[hidelinks]{hyperref}
\makeatletter
\def\ps@pprintTitle{\let\@oddhead\@empty\let\@evenhead\@empty\let\@oddfoot\@empty\let\@evenfoot\@oddfoot}
\makeatother

\begin{document}

\begin{frontmatter}

\title{Hybrid CutFEM-IFEM model for embedded and immersed interfaces with triply-cut cells for the 2D Poisson problem}

\author[lsu-ceds]{Darsh Nathawani}

\author[lsu-cee]{Christopher Kees\corref{cor1}}
\ead{cekees@lsu.edu}
\cortext[cor1]{Corresponding author}

\affiliation[lsu-ceds]{organization={Coastal Ecosystem Design Studio, Louisiana State University},
            city={Baton Rouge},
            postcode={70803},
            state={LA},
            country={USA}}

\affiliation[lsu-cee]{organization={Department of Civil and Environmental Engineering, Louisiana State University},
            city={Baton Rouge},
            postcode={70803},
            state={LA},
            country={USA}}

\begin{abstract}
     We present a hybrid CutFEM--IFEM method for 2D Poisson problems with an immersed fluid--fluid interface and an embedded fluid--solid interface on a single unfitted mesh. The proposed method leverages the strengths of both unfitted finite element methods and immersed interface techniques to handle complex interface geometries and strong discontinuities efficiently. A composite equivalent polynomial representation of the Heaviside and Dirac delta functions is used to integrate on triply-cut cells. Moreover, the IFEM basis functions are enriched with the lift functions exactly representing solution and flux jumps. Numerical examples show optimal convergence for linear and quadratic elements, with accuracy comparable to that of a mesh fitted to both interfaces.
\end{abstract}


\end{frontmatter}

\section{Introduction}

Many coastal and environmental engineering problems involve interactions between multiple fluid phases and solid structures, leading to complex interface dynamics such as breaking waves over beaches~\cite{dimasSurfaceRollerModel2009}, flows through porous and vegetated media~\cite{woodModelingTurbulentFlows2020}, and air--water flows around coastal structures~\cite{keesSimpleAccurateEfficient2022}. In these applications the air--water interface moves and deforms continuously, while the solid boundary may be geometrically complex or itself in motion. The interfaces are commonly tracked implicitly with level set functions (e.g.,~\cite{osherFrontsPropagatingCurvaturedependent1988,sussmanLevelSetApproach1994,keesConservativeLevelSet2011,keesSimpleAccurateEfficient2022}), and generating a mesh that conforms to every interface at every time step is costly and fragile. This motivated the development of unfitted methods, in which the simulation domain resides on a fixed background mesh and the interfaces are allowed to cut the mesh elements. The different ways of performing the special computations on these cut cells lead to the different unfitted methods discussed below.

Unfitted methods have gained significant attention in recent years due to their ability to handle complex geometries and evolving interfaces without the need for mesh conformity. Peskin's immersed boundary method (IBM), developed to study blood flow in the heart, is one of the earliest and most well-known unfitted methods~\cite{peskinNumericalAnalysisBlood1977}. The IBM spreads interfacial forces onto the grid with regularized delta functions and has since been extended to particulate and fluid--structure interaction problems~\cite{uhlmannImmersedBoundaryMethod2005,kempeImprovedImmersedBoundary2012,griffithImmersedMethodsFluid2020}. LeVeque and Li developed the immersed interface method (IIM) to improve the accuracy of interface representation in finite difference schemes~\cite{levequeImmersedInterfaceMethod1994}. IIM introduces corrections to the finite difference stencils near the interface to account for the discontinuities in the coefficients and solution gradients, thereby achieving higher accuracy compared to standard IBM. Their approach was adopted for Stokes, acoustics and elasticity problems, demonstrating its versatility and effectiveness in handling complex interface conditions~\cite{levequeImmersedInterfaceMethods1997,ZhangImmersedInterfaceMethod1997,WiegmannImmersedInterfaceMethod1998,leeImmersedInterfaceMethod2003}. This idea was later implemented with the finite element method by modifying the basis functions on the elements cut by the interface, leading to the development of the immersed finite element method (IFEM)~\cite{liNewCartesianGrid2003}. Since IFEM incorporates the flux discontinuity exactly, the method preserves the sharp jump instead of smearing it across the interface, while keeping the number and location of the degrees of freedom identical to those of the background mesh. This makes IFEM particularly suitable for problems with strong discontinuities and complex interface geometries, providing both accuracy and computational efficiency. The IFEM basis functions are, however, discontinuous across element edges cut by the interface. This non-conformity motivated the symmetric and consistent immersed FEM (SCIFEM)~\cite{jiSymmetricConsistentImmersed2014} and the partially penalized immersed FEM (PPIFEM)~\cite{LinPartiallyPenalizedImmersed2015}, which add edge terms to restore consistency, as well as augmented formulations for non-homogeneous jumps~\cite{jiNewAugmentedImmersed2016}. Higher-order IFEM formulations have also been developed to further improve the accuracy of the method~\cite{AdjeridHigherDegreeImmersed2014,AdjeridHigherDegreeImmersed2018,GuoHigherDegreeImmersed2019,GuzmanHigherorderFiniteElement2016,AdjeridHighOrderGeometry2024}.

A parallel line of work leads to the cut finite element method (CutFEM). Hansbo and Hansbo proposed an unfitted method for elliptic interface problems in which the solution is represented by separate finite element functions on each side of the interface and the interface conditions are imposed weakly with Nitsche's method~\cite{hansboUnfittedFiniteElement2002}. The same idea applied to an embedded boundary yields fictitious domain methods in which Dirichlet conditions are imposed weakly on a boundary that cuts the mesh~\cite{burmanFictitiousDomainFinite2012}. We refer readers to~\cite{burmanCutFEMDiscretizingGeometry2015,burmanCutFiniteElement2019,burmanCutFiniteElementMethods2025} for a detailed discussion of the developments in CutFEM. For an embedded solid boundary, Nitsche's method is a natural choice because the solid region is simply removed from the domain and only a boundary condition must be imposed.

One of the biggest technical challenges in all unfitted methods is computing integrals over the parts of a cut cell that lie inside the physical domain and over the interface segments within the cell. The most common approach subdivides each cut cell into conforming sub-triangles and applies standard quadrature on each of them~\cite{moesExtendedFiniteElement2017,badiaGeometricalDiscretisationsUnfitted2022}. Alternatively, moment fitting constructs cell-specific quadrature weights~\cite{sudhakarQuadratureSchemesArbitrary2013}. Quadrature rules for generalized functions integrate Heaviside and Dirac delta integrands directly~\cite{holdychQuadratureRulesTriangular2008,tornbergMultiDimensionalQuadratureSingular2002}. Dimension-reduction and height-function approaches handle implicitly defined domains~\cite{olshanskiiNumericalIntegrationImplicitly2016,sayeHighorderQuadratureMulticomponent2022}, and closed-form expressions are available for planar cuts~\cite{aulisaExactSubdomainEmbedded2023}. The equivalent polynomial approach replaces the discontinuous Heaviside function by a polynomial that yields the same integrals against all polynomials of a given degree, so that the standard element quadrature can be reused unchanged~\cite{venturaEliminationQuadratureSubcells2006,venturaEquivalentPolynomialsQuadrature2015,abedianEquivalentLegendrePolynomials2019}. Kees et al.\ showed that equivalent polynomials lead to simple, accurate, and efficient embedded finite element methods for fluid--solid interaction with a single level set~\cite{keesSimpleAccurateEfficient2022}.

Problems with three or more phases introduce cells cut by more than one interface and triple points at which the interfaces meet. In the XFEM setting, Zlotnik and D\'{\i}ez used an ordered hierarchy of level sets to describe $n$-phase flows without overlaps or voids~\cite{zlotnikHierarchicalXFEMNphase2009}. Multimesh finite element methods treat multiple overlapping meshes with Nitsche coupling and must integrate over cells intersected by several mesh boundaries~\cite{johanssonMultimeshFiniteElement2019}. For elliptic interface problems with triple junction points, Wang et al.\ proposed a numerical method for the three-dimensional case~\cite{wangNumericalMethodSolving2018}, and Chen et al.\ developed IFEM and partially penalized IFEM formulations for multiple subdomains in two dimensions~\cite{chenImmersedFiniteElement2019,chenBilinearPartiallyPenalized2020}. In these works all interfaces separate materials of the same kind and carry jump conditions, and none of the subdomains is removed from the computational domain. Saye's quadrature handles junctions of multiple polynomial level sets to high order~\cite{sayeHighorderQuadratureMulticomponent2022}. However, the equivalent polynomial technique has so far been developed only for cells cut by a single interface. To the best of our knowledge, no prior work couples an IFEM treatment of an immersed interface with a Nitsche treatment of an embedded boundary within the same cut cell.

In this work, we focus on a three-phase elliptic interface problem, with two fluid phases and one solid phase. We distinguish two kinds of interfaces: an \emph{immersed} interface separates two fluids with different material properties, across which the solution and its flux satisfy jump conditions, and an \emph{embedded} interface is the boundary of a solid region that is removed from the computational domain, on which a boundary condition is prescribed. Both kinds frequently meet inside the same background cell, for instance where the free surface touches a structure. The primary objective is to accurately capture the interface dynamics by combining the strengths of IFEM on the fluid--fluid interface with a Nitsche-based CutFEM formulation on the fluid--solid interface, ensuring both accuracy and computational efficiency in the presence of strong discontinuities. The main contributions are:
\begin{itemize}
    \item a hybrid CutFEM--IFEM formulation, for linear and quadratic elements, that strongly enforces the fluid--fluid jump conditions through the IFEM basis and weakly enforces the fluid--solid boundary condition through Nitsche's method, without adding degrees of freedom;
    \item an extension of equivalent polynomial integration to triply-cut cells;
    \item a consistency correction for triply-cut cells, which accounts for the boundary flux on the element edges removed by the embedded interface and restores optimal convergence;
    \item a treatment of non-homogeneous solution and flux jumps by adding solution-jump and flux-jump functions to the trial function, so that the exact solution remains representable and the error does not depend on the size of the jumps;
    \item a test suite that covers Dirichlet and Robin conditions on the embedded interface, and coefficient contrasts up to $1000$, providing a numerical verification of optimal convergence of order $p+1$ in the $L^2$ and $L^\infty$ norms for $p=1,2$, for homogeneous and non-homogeneous solution and flux jumps, and a comparison with a mesh fitted to both interfaces; with the Robin condition, the quadratic rate is limited to two by the linear representation of the embedded boundary.
\end{itemize}

The rest of the paper is organized as follows. The model problem and its level set description are stated in Section~\ref{sec:model}. Section~\ref{sec:discretization} presents the hybrid discretization, including the IFEM basis construction and the Nitsche formulation. Section~\ref{sec:integration} describes the integration on regular and triply-cut cells. Section~\ref{sec:results} reports the numerical results, and Section~\ref{sec:conclusion} summarizes the conclusions.

\section{Model Problem}
\label{sec:model}

\begin{figure}[htbp]
    \centering
    \includegraphics[width=0.55\linewidth]{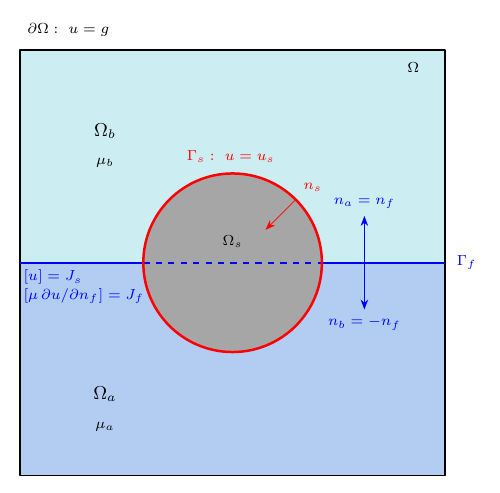}
    \caption{Schematic of the model problem given in Eq.~\eqref{eq:problem_equation}.}
    \label{fig:domain-schematic}
\end{figure}

Consider an elliptic problem in a 2D bounded domain $\Omega$ with a fluid--fluid interface represented by $\Gamma_f$ and a fluid--solid interface represented by $\Gamma_s$. The domain is a union of subdomains $\Omega_f$ and $\Omega_s$, corresponding to the fluid and solid regions, respectively. Furthermore, the fluid domain $\Omega_f$ is made of two fluid regions $\Omega_a$ and $\Omega_b$ separated by the interface $\Gamma_f$. The schematic of the domain and the interfaces is shown in Fig.~\ref{fig:domain-schematic}. The governing equation along with the boundary and interface conditions for the problem is given by
\begin{equation}
    \label{eq:problem_equation}
    \begin{aligned}
        -\nabla \cdot (\mu \nabla u) &= f \quad \text{in } \Omega_f, \\
        u &= g \quad \text{on } \partial\Omega, \\
        u &= u_s \quad \text{on } \Gamma_s, \\
        [u] &= J_s \quad \text{on } \Gamma_f, \\
        \left[\mu \frac{\partial u}{\partial n_f}\right] &= J_f \quad \text{on } \Gamma_f.
    \end{aligned}
\end{equation}
Here, $\mu$ is discontinuous across the interface $\Gamma_f$. We assume that $\mu$ is a piecewise positive constant function, taking values $\mu_a$ and $\mu_b$ in the subdomains $\Omega_a$ and $\Omega_b$, respectively.

The jump conditions across the fluid--fluid interface $\Gamma_f$ are represented by $[u]$ and $[\mu \frac{\partial u}{\partial n_f}]$, which denote the jump in the solution and the jump in the flux, respectively, where $[w] = w_b - w_a$ and $n_f$ is the unit normal on $\Gamma_f$ pointing from $\Omega_a$ into $\Omega_b$. On $\Gamma_s$, $n_s$ denotes the unit normal pointing out of $\Omega_f$, into the solid. The functions $J_s$ and $J_f$ are the prescribed solution and flux jumps. The function $f$ represents the source term, and $g$ is the Dirichlet boundary condition on the outer boundary $\partial\Omega$. The function $u_s$ represents the prescribed solution on the fluid--solid interface $\Gamma_s$.

The interfaces are represented by level set functions, defined as the signed distance from the interface. Since we have two interfaces $\Gamma_f$ and $\Gamma_s$, we use two level set functions, $\Phi_f$ and $\Phi_s$, respectively. The level set function $\Phi_f$ is positive in the subdomain $\Omega_a$ and negative in the subdomain $\Omega_b$, so that $n_f = -\nabla\Phi_f/|\nabla\Phi_f|$ on $\Gamma_f$. The interface $\Gamma_f$ is the zero level set of $\Phi_f$. Similarly, the level set function $\Phi_s$ is positive in the subdomain $\Omega_f$ and negative in the subdomain $\Omega_s$, and $\Gamma_s$ is the zero level set of $\Phi_s$, so that the unit normal $n_s$ pointing out of $\Omega_f$ is $n_s = -\nabla\Phi_s/|\nabla\Phi_s|$ on $\Gamma_s$. The level set functions are used to identify the cut cells in the triangulation $\mathcal{T}_h$ of the domain $\Omega$.

\section{Discretization} \label{sec:discretization}

\subsection{Mesh and cell classification}
\label{sec:mesh}

We seek the approximate solution to Eq.~\eqref{eq:problem_equation} using piecewise polynomials defined on a triangulation $\mathcal{T}_h$ of $\Omega$ that is independent of the interfaces. In the unfitted-FEM terminology, this triangulation is called a background mesh, and the elements that are cut by the interfaces are called cut cells. The discrete interfaces $\Gamma_{f,h}$ and $\Gamma_{s,h}$ are the zero level sets of the piecewise linear interpolants of $\Phi_f$ and $\Phi_s$, so that each interface is a straight segment within a cut cell.

\begin{figure}[t]
    \centering
    \includegraphics[width=0.65\linewidth]{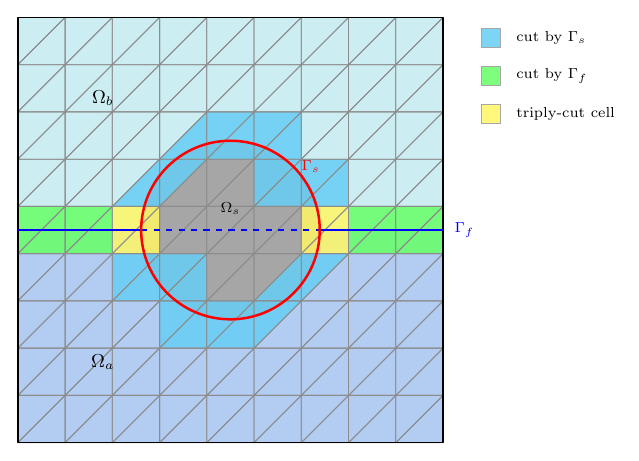}
    \caption{Triangulation of the domain with doubly- and triply-cut cells.}
    \label{fig:triangulation-schematic}
\end{figure}

Only the cells that contain fluid enter the computation. We call them active cells,
\begin{align*}
    \mathcal{T}_h^{\mathrm{act}} = \{ K \in \mathcal{T}_h : |K \cap \Omega_f| > 0 \},
\end{align*}
and classify them as
\begin{align*}
    \mathcal{T}_h^{s} &= \{ K \in \mathcal{T}_h^{\mathrm{act}} : K \cap \Gamma_{s,h} \neq \emptyset,\ K \cap \Gamma_{f,h} = \emptyset \} && \text{(cut by } \Gamma_s \text{ only)}, \\
    \mathcal{T}_h^{f} &= \{ K \in \mathcal{T}_h^{\mathrm{act}} : K \cap \Gamma_{f,h} \neq \emptyset,\ K \cap \Gamma_{s,h} = \emptyset \} && \text{(cut by } \Gamma_f \text{ only)}, \\
    \mathcal{T}_h^{fs} &= \{ K \in \mathcal{T}_h^{\mathrm{act}} : K \cap \Gamma_{f,h} \neq \emptyset,\ K \cap \Gamma_{s,h} \neq \emptyset \} && \text{(cut by both, triply-cut)}, \\
    \mathcal{T}_h^{0} &= \mathcal{T}_h^{\mathrm{act}} \setminus \left( \mathcal{T}_h^{s} \cup \mathcal{T}_h^{f} \cup \mathcal{T}_h^{fs} \right) && \text{(regular)}.
\end{align*}
Cells entirely inside the solid are not active and are discarded. $\mathcal{N}_h$ denotes the set of nodes of $\mathcal{T}_h^{\mathrm{act}}$, and $h_K$ the diameter of $K$. An example mesh with the cut cells is shown in Fig.~\ref{fig:triangulation-schematic}. Here, the gray elements inside $\Omega_s$ are inactive. All other elements are active elements.

\subsection{Function spaces}
\label{sec:ife-space}

We begin with the construction based on local linear functions. On the cells that are not cut by $\Gamma_f$, $K \in \mathcal{T}_h^{0} \cup \mathcal{T}_h^{s}$, we use the standard linear space $V^{ifem}_h(K) = P^1(K)$. On the cells cut by $\Gamma_f$, $K \in \mathcal{T}_h^{f} \cup \mathcal{T}_h^{fs}$, we use the IFEM space $V^{ifem}_h(K)$, whose basis functions are constructed such that they include the flux discontinuity across $\Gamma_{f,h}$~\cite{jiSymmetricConsistentImmersed2014}.

Let $K_a = \{ \mathbf{x} \in K : \Phi_f(\mathbf{x}) > 0 \}$ and $K_b = \{ \mathbf{x} \in K : \Phi_f(\mathbf{x}) < 0 \}$. Then the IFEM space $V^{ifem}_h(K) = \mathrm{span}\{\phi_1, \phi_2, \phi_3\}$, where
\begin{align}
    \phi_i(\mathbf{x}) = \begin{cases}
        \phi_i^a(\mathbf{x}) = a_1 + a_2 x + a_3 y & \text{if } \mathbf{x} \in K_a, \\
        \phi_i^b(\mathbf{x}) = b_1 + b_2 x + b_3 y & \text{if } \mathbf{x} \in K_b,
    \end{cases}
\end{align}
subject to the following constraints:
\begin{equation}
    \label{eq:ifem-basis-constraints}
    \begin{aligned}
        \phi_i(\mathbf{x}_j) &= \delta_{ij}, && j=1,\dots,3, \\
        \phi_i^a &= \phi_i^b && \text{at the two end points of } \Gamma_{f,h}\cap K, \\
        \mu_a \frac{\partial \phi_i^a}{\partial n_f} &= \mu_b \frac{\partial \phi_i^b}{\partial n_f} && \text{at the midpoint of } \Gamma_{f,h}\cap K.
    \end{aligned}
\end{equation}

\begin{figure}[t]
    \centering
    \includegraphics[width=\linewidth]{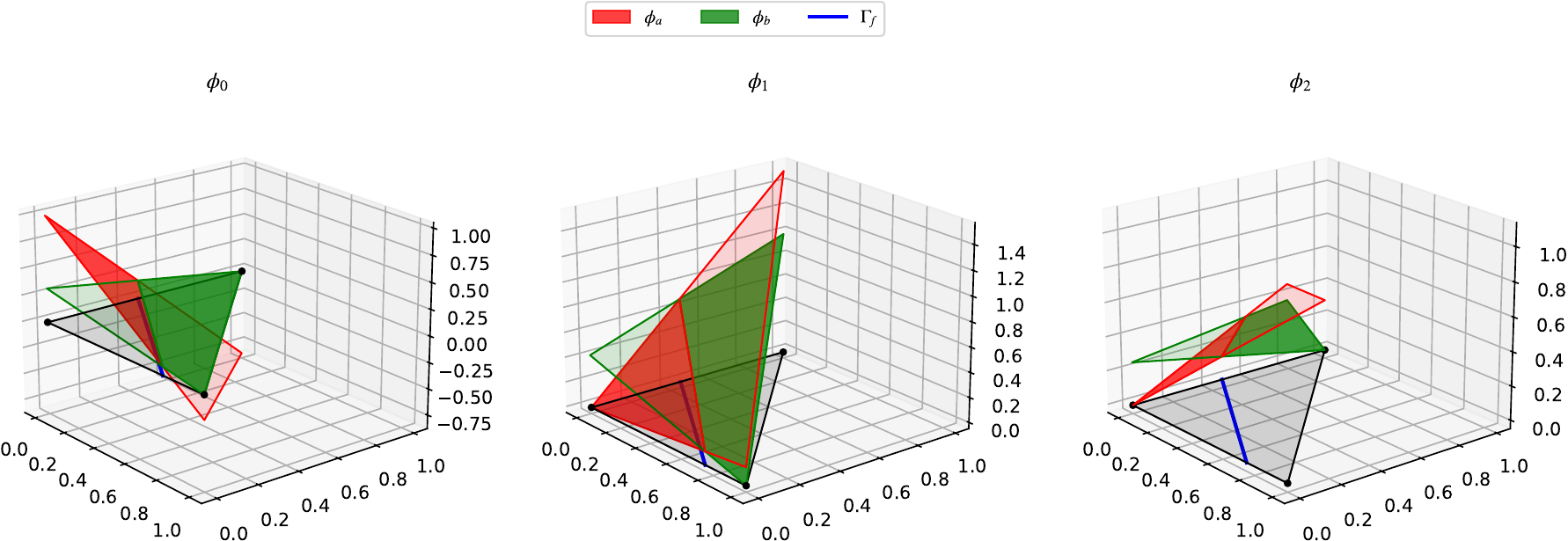}
    \caption{Linear basis functions constructed on a cut cell.}
    \label{fig:linear-basis-on-cut-cell}
\end{figure}

The first constraint ensures the Kronecker delta property of the basis functions, the second constraint enforces the continuity of the basis functions across the interface $\Gamma_{f,h}$, and the third constraint enforces the continuity of the flux across the interface $\Gamma_{f,h}$. Collectively, we have a total of six equations to solve for the six unknowns $a_1, a_2, a_3, b_1, b_2, b_3$. An example of the basis functions constructed on a cut cell is shown in Fig.~\ref{fig:linear-basis-on-cut-cell}. We define the function spaces
\begin{align}
    V_h &= \left\{ v^h : v^h|_K \in V^{ifem}_h(K)\ \forall K \in \mathcal{T}_h^{\mathrm{act}},\ v^h = g \text{ on } \partial\Omega \right\}, \\
    V_h^0 &= \left\{ v^h : v^h|_K \in V^{ifem}_h(K)\ \forall K \in \mathcal{T}_h^{\mathrm{act}},\ v^h = 0 \text{ on } \partial\Omega \right\},
\end{align}
for the trial and test functions, respectively; in both, $g$ and $0$ are imposed at the nodes on $\partial\Omega$.

Note that the IFEM space $V^{ifem}_h(K)$ reduces to the standard linear space when there is no jump in the coefficients, i.e.\ $\mu_a = \mu_b$. However, the IFEM basis functions are discontinuous across the edges cut by $\Gamma_f$, hence $V_h$ is non-conforming. On a triply-cut cell $K \in \mathcal{T}_h^{fs}$, the IFEM basis functions are constructed in the same way, using $\Phi_f$, unaware of $\Gamma_s$. $\Gamma_s$ enters only through the weak form explained in Section~\ref{sec:weak-form}.

\vspace{5mm}
\noindent
\textit{Quadratic elements}
\vspace{5mm}

\begin{figure}[!b]
    \centering
    \includegraphics[width=\linewidth]{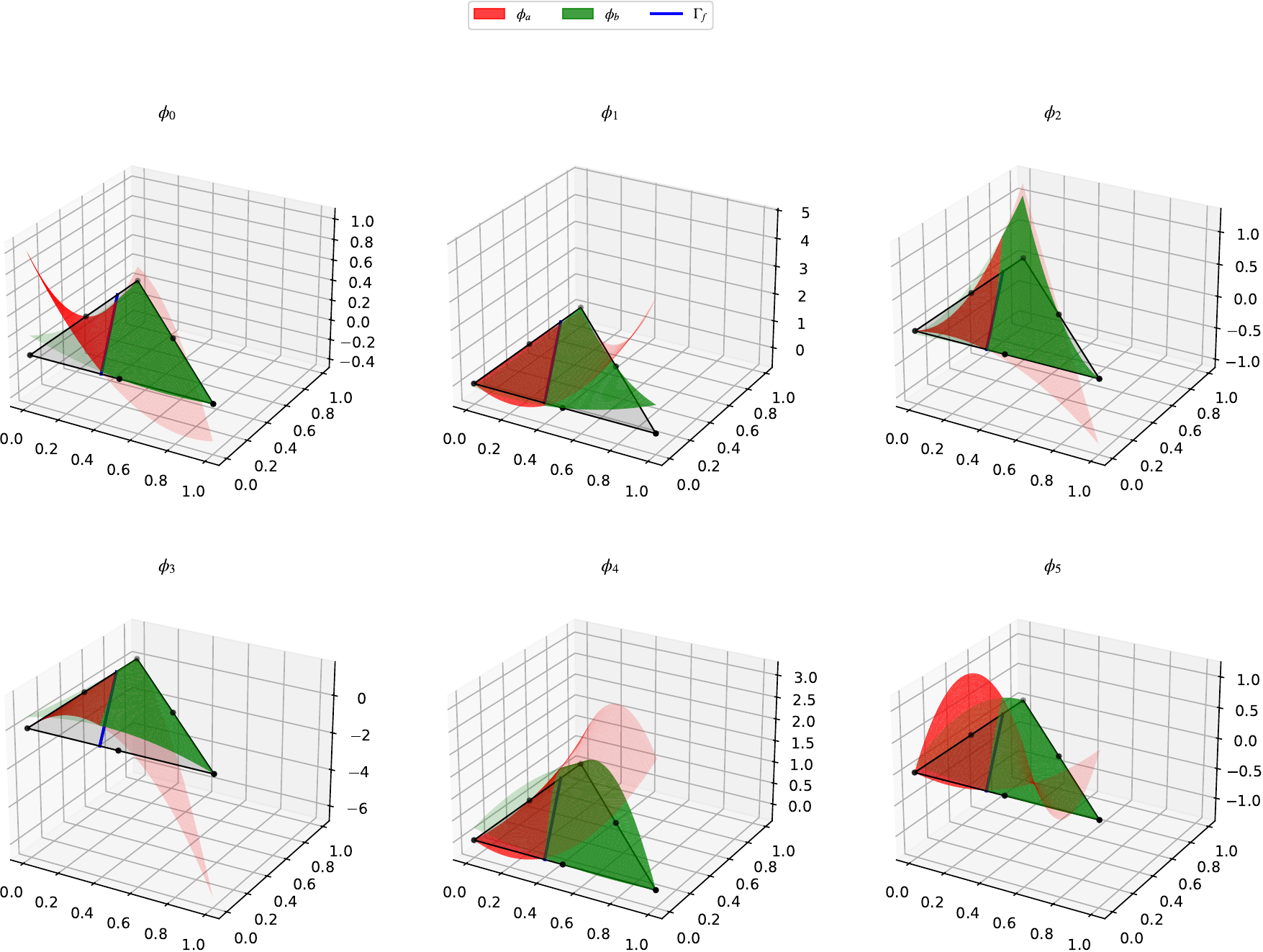}
    \caption{Quadratic basis functions constructed on a cut cell.}
    \label{fig:quadratic-basis-on-cut-cell}
\end{figure}

The construction extends to polynomials of degree $p=2$. On a cell not cut by $\Gamma_f$ the space is $P^2(K)$. On a cell cut by $\Gamma_{f,h}$, each basis function is a pair of quadratics,
\begin{align}
    \phi_i(\mathbf{x}) = \begin{cases}
        \phi_i^a(\mathbf{x}) = a_1 + a_2 x + a_3 y + a_4 xy + a_5 x^2 + a_6 y^2 & \text{in } K_a, \\
        \phi_i^b(\mathbf{x}) = b_1 + b_2 x + b_3 y + b_4 xy + b_5 x^2 + b_6 y^2 & \text{in } K_b,
    \end{cases}
\end{align}
where $i=1,\ldots,6$. The twelve coefficients are fixed by twelve conditions:
\begin{equation}
    \label{eq:p2-constraints}
    \begin{aligned}
        \phi_i(\mathbf{x}_j) &= \delta_{ij}, && j=1,\dots,6, \\
        \phi_i^a &= \phi_i^b && \text{at the two end points and the midpoint of } \Gamma_{f,h}\cap K, \\
        \mu_a \frac{\partial \phi_i^a}{\partial n_f} &= \mu_b \frac{\partial \phi_i^b}{\partial n_f} && \text{at the two end points of } \Gamma_{f,h}\cap K, \\
        \mu_a \frac{\partial^2 \phi_i^a}{\partial n_f^2} &= \mu_b \frac{\partial^2 \phi_i^b}{\partial n_f^2} && \text{at the midpoint of } \Gamma_{f,h}\cap K .
    \end{aligned}
\end{equation}
The last condition can be derived from the PDE consistency. Decomposing the Laplacian on both sides, and dropping the curvature term and the tangential term due to linear cut and continuity of $\phi$, leads to the normal derivative term given above. The six basis functions of a cut cell are shown in Fig.~\ref{fig:quadratic-basis-on-cut-cell}.

\vspace{5mm}
\noindent
\textit{Solution jump}
\vspace{5mm}

\begin{figure}[!ht]
    \centering
    \begin{subfigure}{0.45\linewidth}
        \centering
        \includegraphics[width=\linewidth]{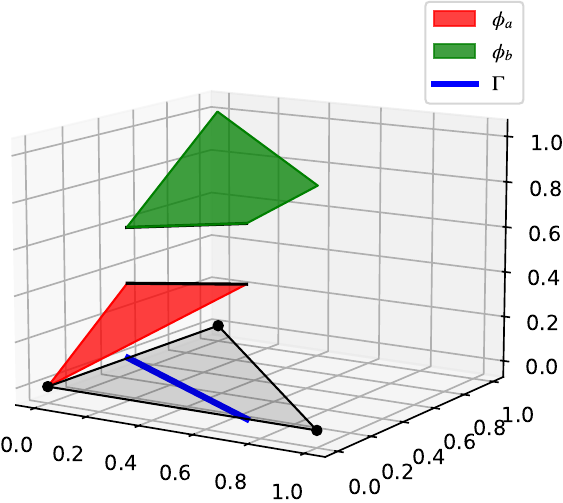}
        \caption{Solution jump, $p = 1$.}
        \label{fig:jump-lift-p1}
    \end{subfigure}
    \hfill
    \begin{subfigure}{0.45\linewidth}
        \centering
        \includegraphics[width=\linewidth]{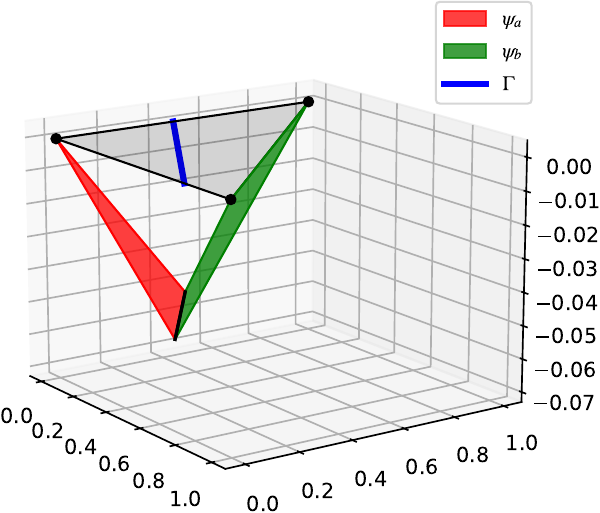}
        \caption{Flux jump, $p = 1$.}
        \label{fig:flux-lift-p1}
    \end{subfigure}
    \\[2ex]
    \begin{subfigure}{0.45\linewidth}
        \centering
        \includegraphics[width=\linewidth]{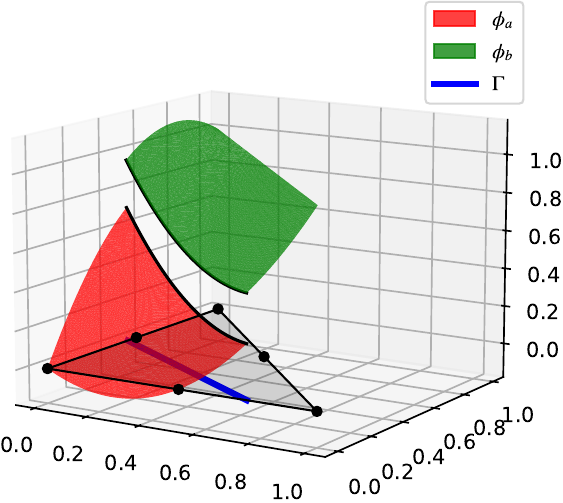}
        \caption{Solution jump, $p = 2$.}
        \label{fig:jump-lift-p2}
    \end{subfigure}
    \hfill
    \begin{subfigure}{0.45\linewidth}
        \centering
        \includegraphics[width=\linewidth]{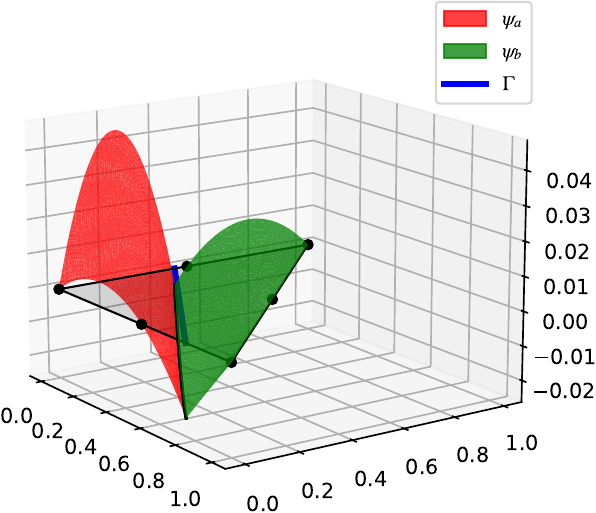}
        \caption{Flux jump, $p=2$.}
        \label{fig:flux-lift-p2}
    \end{subfigure}
    \caption{The solution and flux jump enrichment functions on a cut cell, for $p=1$ (top) and $p=2$ (bottom). $\mu_a=1$, $\mu_b=5$, and $J_s = 0.25$.}
    \label{fig:lifts}
\end{figure}

The functions in $V_h$ are continuous across $\Gamma_{f,h}$. The prescribed solution jump $[u] = J_s$ is handled by augmenting the trial function with a jump function, $u^h + j_h$ with $u^h \in V_h$. On $K \in \mathcal{T}_h^{f} \cup \mathcal{T}_h^{fs}$,
\begin{align}
    \label{eq:jump-lifting}
    j_h = \begin{cases}
        \sum_{k} j_k^a\, \phi_k^a & \text{in } K_a, \\
        \sum_{k} j_k^b\, \phi_k^b & \text{in } K_b,
    \end{cases}
    \qquad
    (j_k^a, j_k^b) = \begin{cases}
        \big(0,\ J_s(\mathbf{x}_k)\big) & \text{if } \mathbf{x}_k \in \overline{K}_a, \\
        \big(-J_s(\mathbf{x}_k),\ 0\big) & \text{if } \mathbf{x}_k \in \overline{K}_b,
    \end{cases}
\end{align}
and $j_h = 0$ on all other cells. Fig.~\ref{fig:jump-lift-p1} and Fig.~\ref{fig:jump-lift-p2} show linear and quadratic functions augmented in this way, respectively, such that the restrictions to $K_a$ and $K_b$ are discontinuous across $\Gamma_{f,h}$ by the prescribed jump, while $u^h$ itself is unchanged.

\vspace{5mm}
\noindent
\textit{Flux jump}
\vspace{5mm}

The IFEM basis functions satisfy the homogeneous flux condition of Eq.~\eqref{eq:ifem-basis-constraints}, as they must for $V_h$ to be a linear space. To represent the non-homogeneous flux jump, the IFEM space needs enrichment functions having the prescribed flux jump across the interface, as explored by Adjerid et al.~\cite{adjeridEnrichedImmersedFinite2023}. The prescribed flux jump enters the weak form as the load term derived in Section~\ref{sec:weak-form}. On $K \in \mathcal{T}_h^{f} \cup \mathcal{T}_h^{fs}$, let $\psi_K = (\psi_K^a, \psi_K^b)$, with $\psi_K^a,\psi_K^b \in P^p(K)$, satisfy the conditions of the basis with zero nodal values and a unit flux jump:
\begin{equation}
    \label{eq:psi}
    \begin{aligned}
        &\psi_K(\mathbf{x}_j) = 0 \ \ \text{for all nodes } \mathbf{x}_j \text{ of } K, \qquad
        \psi_K^a = \psi_K^b \ \ \text{on } \Gamma_{f,h}\cap K, \\[2ex]
        &\mu_b \frac{\partial \psi_K^b}{\partial n_f} - \mu_a \frac{\partial \psi_K^a}{\partial n_f} = 1 \ \ \text{on } \Gamma_{f,h}\cap K, \qquad
        \mu_b \frac{\partial^2 \psi_K^b}{\partial n_f^2} - \mu_a \frac{\partial^2 \psi_K^a}{\partial n_f^2} = 0 \ \ (p=2),
    \end{aligned}
\end{equation}
where the continuity condition is imposed at $p+1$ points of the segment $\Gamma_{f,h}\cap K$ and the flux condition at $p$ points: the midpoint for $p=1$ and the two end points for $p=2$. Since the enrichment functions are zero at the nodes, the solution degrees of freedom are unaffected. The flux-lift function is then represented by
\begin{equation}
    q_h = J_f(\mathbf{x}_k) \, \psi_K  \qquad \mathbf{x}_k \text{ the midpoint of } \Gamma_{f,h}\cap K,
\end{equation}
and $q_h = 0$ on all other cells. Fig.~\ref{fig:flux-lift-p1} and Fig.~\ref{fig:flux-lift-p2} show the linear and quadratic flux enrichment functions, respectively.

With the solution jump and the flux jump enrichment, the trial function now becomes $u^h + j_h + q_h \text{ for } u^h \in V_h$. For $J_s$ and $J_f$ constant on $K$, the augmented trial function has exactly the prescribed jumps, $[u^h + j_h + q_h]=J_s$ and $[\mu\,(\partial u^h/\partial n_f + \partial j_h/\partial n_f +\partial q_h/\partial n_f)]=J_f$ on $\Gamma_{f,h}\cap K$.

\subsection{Weak form} \label{sec:weak-form}

The weak form of the problem is derived by multiplying Eq.~\eqref{eq:problem_equation} by a test function $v^h \in V^0_h$ and integrating over the fluid domain $\Omega_f$. Integrating by parts on $\Omega_a$ and $\Omega_b$ separately gives
\begin{align}
    \label{eq:bilinear-form-split}
    \int_{\Omega_f} -\nabla\cdot(\mu \nabla u)\, v^h \, dx
    &= \sum_{i \in \{a,b\}} \int_{\Omega_i} \mu_i \nabla u_i \cdot \nabla v^h \, dx
    - \sum_{i \in \{a,b\}} \int_{\partial\Omega_i} \mu_i \frac{\partial u_i}{\partial n_i}\, v^h \, ds,
\end{align}
where $u_i = u|_{\Omega_i}$ and $n_i$ is the outward unit normal of $\Omega_i$. The boundary $\partial\Omega_i$ consists of parts of $\partial\Omega$, $\Gamma_s$ and $\Gamma_f$, which we treat in turn.
\begin{itemize}
    \item On $\partial\Omega$, the integral vanishes because $v^h = 0$ there.
    \item On $\Gamma_f$, with $n_a = n_f$ and $n_b = -n_f$, the two integrals combine into a flux jump term given by
    \begin{align*}
        -\int_{\Gamma_f} \mu_a \frac{\partial u_a}{\partial n_a}\, v^h \, ds - \int_{\Gamma_f} \mu_b \frac{\partial u_b}{\partial n_b}\, v^h \, ds &= \int_{\Gamma_f} \left[ \mu \frac{\partial u}{\partial n_f} \right] v^h \, ds \\
        &= \int_{\Gamma_f} J_f\, v^h \, ds,
    \end{align*}
    which is known data and moves to the right-hand side.
    \item On $\Gamma_s$, the outward normal of $\Omega_f$ is $n_s$, and the boundary integral is $-\int_{\Gamma_s} \mu \frac{\partial u}{\partial n_s} v^h\, ds$. Nitsche's method is used to weakly enforce the Dirichlet boundary condition $u = u_s$~\cite{hansboUnfittedFiniteElement2002} by keeping this integral, adding a symmetric term that vanishes for $u = u_s$, and adding a penalty term to restore coercivity:
    \begin{align*}
        \int_{\Gamma_{s}} \left( -\mu \frac{\partial u^h}{\partial n_s} v^h - \mu \frac{\partial v^h}{\partial n_s} (u^h - u_s) + \frac{\gamma_N \mu}{h_K} (u^h - u_s) v^h \right) ds .
    \end{align*}
    If instead a natural boundary condition is given, the boundary flux term is replaced by the given data, eliminating the need for Nitsche's terms.
\end{itemize}

For notational simplicity, we write $u^h$ for the trial function $u^h + j_h + q_h$ in the forms below; since $j_h$ and $q_h$ are known, their contributions move to the right-hand side. The discrete problem reads: find $u^h \in V_h$ such that
\begin{align}
    \label{eq:hybrid-discrete}
    a(u^h, v^h) + c(u^h, v^h) + s(u^h,v^h) + s_I(u^h,v^h) = l(v^h) \qquad \forall v^h \in V^0_h,
\end{align}
where
\begin{equation}
\label{eq:a-and-l}
\begin{aligned}
    a(u^h,v^h) &= \int_{\Omega_f} \mu \nabla u^h \cdot \nabla v^h \, dx
    + \int_{\Gamma_{s}} \left( -\mu \frac{\partial u^h}{\partial n_s} v^h - \mu \frac{\partial v^h}{\partial n_s} u^h + \frac{\gamma_N \mu}{h_K} u^h v^h \right) ds,  \\[2ex]
    l(v^h) &= \int_{\Omega_f} f v^h \, dx
    - \int_{\Gamma_{f} \cap \Omega_f} J_f\, v^h \, ds
    + \int_{\Gamma_{s}} \left( -\mu \frac{\partial v^h}{\partial n_s} u_s + \frac{\gamma_N \mu}{h_K} u_s v^h \right) ds,
\end{aligned}
\end{equation}
and $\gamma_N > 0$ is the Nitsche penalty parameter.

To provide additional stability and bound the condition number independently of the cut location, CutFEM adds a \emph{ghost penalty} term $s(u^h, v^h)$. It penalizes the jump of the normal derivative across the faces of the cells that $\Gamma_s$ cuts, and thereby ties the gradient in a sliver cell to the gradient in its neighbors. Let $\mathcal{E}^s$ be the set of all $\Gamma_s$-cut edges and $h_e$ the length of the edge $e \in \mathcal{E}^s$. The ghost penalty term is
\begin{equation}
    \label{eq:ghost-penalty}
    s(u^h, v^h) = \gamma_g \sum_{e\in\mathcal{E}^{s}} h_e \int_e \left[\!\left[ \frac{\partial u^h}{\partial n_e} \right]\!\right]_e \left[\!\left[ \frac{\partial v^h}{\partial n_e} \right]\!\right]_e ds ,
\end{equation}
with the edge jump $[\![\cdot]\!]_e$ defined in Eq.~\eqref{eq:avg-jump}.

The terms $c(u^h,v^h)$ and $s_I(u^h,v^h)$ account for the discontinuity of the basis functions across the $\Gamma_f$-cut edges; they are defined in the following section.

\subsection{Weak form consistency in \texorpdfstring{$\Gamma_f$}{Gamma\_f}-cut cells}
\label{sec:consistency}

In the standard finite element method, the weak form is obtained by multiplying the equation by a test function and integrating by parts on each element. Summing over the elements leaves, on every interior edge, a boundary term that couples the flux of the solution with the values of the test function on the two sides of the edge. These terms vanish because the exact flux is continuous across the edge and the test functions are continuous, so that the weak form can be stated over the whole domain without any edge integrals. In IFEM the basis functions are constructed to satisfy the jump conditions on the interface $\Gamma_f$. Hence, the test functions are only piecewise smooth and are not continuous across the edges cut by the interface, so these edge integrals do not vanish and the weak form without them is not consistent. This non-conformity of the basis functions led to the construction of the symmetric and consistent immersed finite element method (SCIFEM) in~\cite{jiSymmetricConsistentImmersed2014} and the partially penalized immersed finite element method (PPIFEM) in~\cite{LinPartiallyPenalizedImmersed2015}, which retain the edge terms. Nevertheless, the standard IFEM still provides optimal convergence because the non-conformity vanishes on the vertices. However, in the case of a triply-cut cell, the fluid--solid interface $\Gamma_s$ chips away a part of the shared edge cut by the fluid--fluid interface $\Gamma_f$. As a result, the boundary integral over the shared cut edge adds a nontrivial contribution to the weak form. Therefore, a standard non-conforming IFEM approach leads to inconsistency.

Let $\mathcal{E}^{f}$ be the set of all edges of $\mathcal{T}_h^{\mathrm{act}}$ that are cut by $\Gamma_{f}$ and have a fluid part, $|e\cap\Omega_f|>0$. Let $n_e$ be a unit normal of the edge $e \in \mathcal{E}^f$, pointing out of the cell under consideration. The contribution from the boundary integrals on these edges is
\begin{equation}
    \label{eq:non-symmetric-ifem-term}
    c(u^h,v^h) = - \sum_{e\in\mathcal{E}^{f}} \int_{e\cap\Omega_f} \mu \frac{\partial u^h}{\partial n_e} \, v^h \, ds .
\end{equation}

\begin{figure}[htbp]
    \centering
    \includegraphics[width=0.3\textwidth]{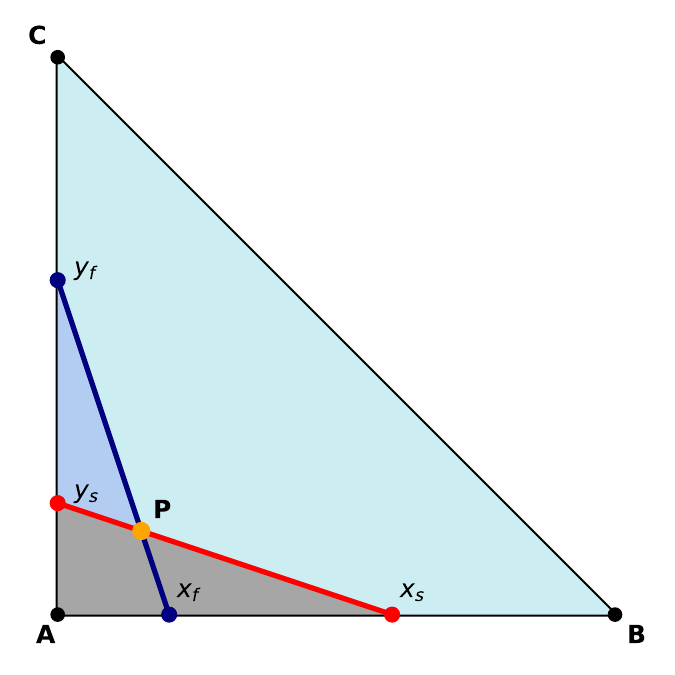}
    \caption{Example triply-cut element.}
    \label{fig:example-triply-cut}
\end{figure}

Alternatively, the consistency term can be treated with symmetric contributions on the interior edges. Let $\mathcal{E}^{f} = \mathcal{E}_{int}^f \cup \mathcal{E}_{ext}^f$, the union of the interior and exterior $\Gamma_f$-cut edges. An edge $e\in\mathcal{E}_{int}^{f}$ is shared by two cells $K_1$ and $K_2$ of $\mathcal{T}_h^{f}\cup\mathcal{T}_h^{fs}$. Let $n_i$ be the outward unit normal of $K_i$ on $e$; we set $n_e=n_1$, so that $n_2=-n_e$. Following Ji et al.~\cite{jiSymmetricConsistentImmersed2014}, the average and the jump across $e$ are
\begin{align}
    \label{eq:avg-jump}
    \{\!\{w\}\!\}_e = \tfrac12\left( w|_{K_1} + w|_{K_2} \right), \qquad
    [\![w]\!]_e = w|_{K_1} - w|_{K_2}.
\end{align}
The double brackets denote the jump across an edge, to distinguish it from the jump $[\cdot]$ across $\Gamma_f$. The edge terms are
\begin{equation}
    \label{eq:scifem-consistency-term}
    \begin{aligned}
        c(u^h,v^h) = &- \sum_{e\in\mathcal{E}_{int}^{f}} \int_{e\cap\Omega_f} \left( \left\{\!\!\left\{ \mu\frac{\partial u^h}{\partial n_e} \right\}\!\!\right\}_e [\![v^h]\!]_e + \left\{\!\!\left\{ \mu\frac{\partial v^h}{\partial n_e} \right\}\!\!\right\}_e [\![u^h]\!]_e \right) ds \\
        &- \sum_{e\in\mathcal{E}_{ext}^{f}} \int_{e\cap\Omega_f} \mu \frac{\partial u^h}{\partial n_e} \, v^h \, ds .
    \end{aligned}
\end{equation}
Note that the edges $e \in \mathcal{E}^f$ that are in triply-cut cells $K \in \mathcal{T}_h^{fs}$ are cut by $\Gamma_s$ as well. Therefore, the integrals above are applied to the fluid part of those edges only. As an example, consider the triply-cut cell shown in Fig.~\ref{fig:example-triply-cut}. The fluid-containing parts are the segment from $x_s$ to $B$ on the edge $AB$, and the segment from $y_s$ to $C$ on the edge $AC$. All the triply-cut cell configurations shown in Fig.~\ref{fig:triply-cut-cell-configs} have similar consistency contributions.

Note that if the fluid--fluid interface $\Gamma_f$ does not intersect the domain boundary $\partial \Omega$, the consistency term given in Eq.~\eqref{eq:scifem-consistency-term} is symmetric. We use Eq.~\eqref{eq:scifem-consistency-term} for all the numerical examples in Section~\ref{sec:results}.

\vspace{5mm}
\noindent
\textit{Interior penalty on the $\Gamma_f$-cut edges}
\vspace{5mm}

In addition to the edge terms of Eq.~\eqref{eq:scifem-consistency-term} to ensure consistency, we add a penalty term on the jumps across the interior $\Gamma_f$-cut edges, as in the partially penalized IFEM~\cite{LinPartiallyPenalizedImmersed2015}. This ensures the coercivity of the discrete form.
\begin{align}
    \label{eq:interior-penalty}
    s_I(u^h,v^h) = \sum_{e\in\mathcal{E}_{int}^{f}} \frac{\gamma_I}{h_e} \sum_{i\in\{a,b\}} \mu_i \int_{e\cap\Omega_i\cap\Omega_f} [\![u_i^h]\!]_e\,[\![v_i^h]\!]_e \, ds ,
\end{align}
where $\gamma_I \ge 0$ is the penalty parameter and $u_i^h$ and $v_i^h$ are the restrictions of $u^h$ and $v^h$ to $\Omega_i$. The penalty acts separately on the parts of each edge in $\Omega_a$ and $\Omega_b$, and each part is weighted by its coefficient $\mu_i$, as in Guo and Lin~\cite{GuoHigherDegreeImmersed2019}, so that the penalty scales with the bulk term on both sides of $\Gamma_f$.

\section{Integration in cut cells}
\label{sec:integration}

In order to use the equivalent polynomial approach, we rewrite the discrete weak form in terms of Heaviside and Dirac delta functions as follows:
\begin{align}
     a(u^h,v^h) &= \int_{\Omega} H(\Phi) \mu \nabla u^h \cdot \nabla v^h \, dx
    - \int_{\Omega} \delta(\Phi_s) \mu \frac{\partial u^h}{\partial n_s}\, v^h \, dx \nonumber \\
    & \quad - \int_{\Omega} \delta(\Phi_s) \mu \frac{\partial v^h}{\partial n_s} u^h \, dx
    + \int_{\Omega} \delta(\Phi_s) \frac{\gamma_N \mu}{h_K} u^h \, v^h \, dx, \\[2ex]
    l(v^h) &= \int_{\Omega} H(\Phi) f v^h \, dx - \int_{\Omega} \delta(\Phi_f)\, J_f\, v^h \, dx \nonumber \\
    & \quad - \int_{\Omega} \delta(\Phi_s) \mu \frac{\partial v^h}{\partial n_s} u_s \, dx + \int_{\Omega} \delta(\Phi_s) \frac{\gamma_N \mu}{h_K} u_s\, v^h \, dx,
\end{align}
where, in the bulk terms, $\Phi$ is $\Phi_s$ when the cell is cut by the fluid--solid interface $\Gamma_s$ and $\Phi_f$ when the cell is cut by the fluid--fluid interface $\Gamma_f$. Here $H$ is the Heaviside function defined as
\begin{align}
    H(z) = \begin{cases}
        1 & \text{if } z > 0, \\
        0 & \text{if } z < 0,
    \end{cases}
\end{align}
and the Dirac delta function is $\delta(z) = H'(z)$.

The integrand next to the Heaviside and Dirac delta functions is often a polynomial of known degree. As described in~\cite{venturaEquivalentPolynomialsQuadrature2015}, the Heaviside and Dirac delta functions are represented by \textit{equivalent polynomials} with at least the same degree $k$ as the integrand. As a result, we get an integrand of degree $2k$, which can be integrated exactly using appropriate quadrature rules. Consider a reference element $K_r$ that is cut by a linear interface $\Gamma = \{\textbf{x} \mid \Phi(\textbf{x}) = 0\}$, with a polynomial integrand $F(\textbf{x})$ of degree $k$ to integrate. The equivalent polynomial $H_{eq}$ of degree $k$ is defined by
\begin{align*}
    \int_{K_r} H_{eq}\, F(\textbf{x}) \, dx = \int_{K_r} H(\Phi(\textbf{x}))\, F(\textbf{x}) \, dx \qquad \text{for all } F \in P^k(K_r).
\end{align*}
Since this identity must hold for every polynomial of degree up to $k$, we seek $H_{eq}$ in $P^k(K_r)$: we write it as a linear combination of a basis $\{\varphi_i\}$ of $P^k(K_r)$ and compute the coefficients by moment fitting,
\begin{align}
    \int_{K_r} \sum_{i=1}^{\dim P^k(K_r)} \hat{H}_i \varphi_i \varphi_j \, dx &= \int_{K_r} H(\Phi(\textbf{x})) \varphi_j \, dx , \quad j = 1,\dots,\dim P^k(K_r), \\
    \therefore \quad \mathbb{M} \hat{\textbf{H}} &= \textbf{b},
\end{align}
where $\mathbb{M}$ is the mass matrix and $\textbf{b}$ is the right-hand side vector. Note that the mass matrix is independent of the interfaces and can be precomputed. The right-hand side vector $\textbf{b}$ can be computed using the signed distance functions $\Phi_f$ and $\Phi_s$ depending on the type of the cut cell. In the following, we describe the procedure to compute the right-hand side vector $\textbf{b}$ for a regular cell cut by $\Gamma_f$ or $\Gamma_s$, and for a triply-cut cell.
\subsection{Regular cut cell}
\begin{figure}[htbp]
    \centering
    \begin{subfigure}{0.45\textwidth}
        \centering
        \includegraphics[width=0.7\linewidth]{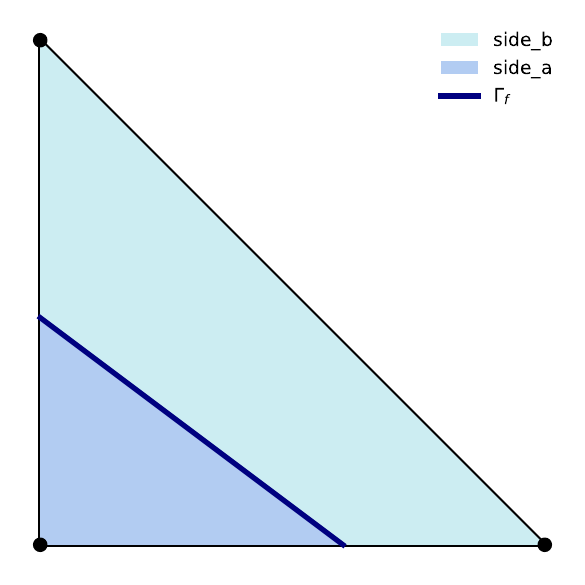}
        \caption{Regular cell cut by $\Gamma_f$.}
        \label{fig:regular-cut-cell-fluid-fluid}
    \end{subfigure}
    \hfill
    \begin{subfigure}{0.45\textwidth}
        \centering
        \includegraphics[width=0.7\linewidth]{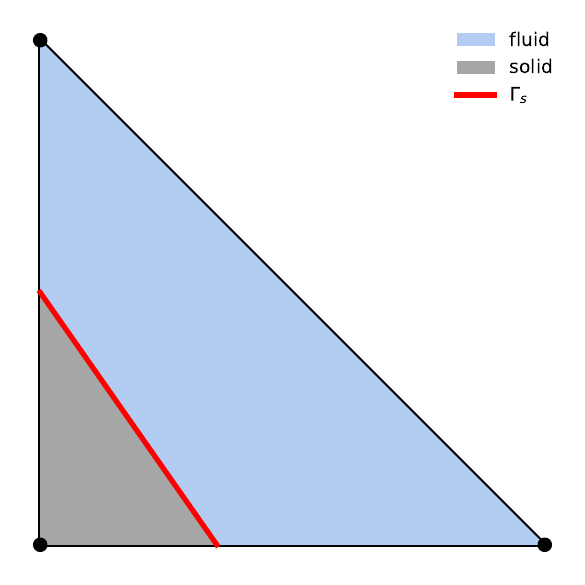}
        \caption{Regular cell cut by $\Gamma_s$.}
        \label{fig:regular-cut-cell-solid}
    \end{subfigure}
    \caption{Example of regular cells cut by a single interface.}
    \label{fig:regular-cut-cell}
\end{figure}

In the case of a cell that is cut by only one interface, we can compute the right-hand side as a function of a single level set function $\Phi_f$ or $\Phi_s$. Without loss of generality, say the interface isolates node $0$ of $K_r$, crossing edges $0$-$1$ and $0$-$2$ at the points $(x_c,0)$ and $(0,y_c)$. Since the level set is a signed distance to a straight interface, each cut point is exactly the zero of the affine interpolant of the level set along that edge:
\begin{align}
    \label{eq:cut-formula}
    x_c = \frac12 - \frac12\frac{\Phi^1+\Phi^0}{\Phi^1-\Phi^0}, \qquad
    y_c = \frac12 - \frac12\frac{\Phi^2+\Phi^0}{\Phi^2-\Phi^0},
\end{align}
where $\Phi^0,\Phi^1,\Phi^2$ are the nodal values of whichever level set cuts the cell. An example of a reference regular cut cell is shown in Fig.~\ref{fig:regular-cut-cell}.

For the fluid--solid cut cell, let $K_s = \{\textbf{x} \in K_r \mid \Phi_s(\textbf{x}) < 0\}$ and let $\{\Phi_s^0, \Phi_s^1, \Phi_s^2\}$ be the nodal values of the level set function $\Phi_s$. Then
\begin{align}
    b_j &= \int_{K_r} (1 - H(\Phi_s)) \varphi_j \, dx = \int_{K_s} \varphi_j \, dx .
\end{align}

Similarly, for the fluid--fluid cut cell, let $K_a = \{\textbf{x} \in K_r \mid \Phi_f(\textbf{x}) > 0\}$ and let $\{\Phi_f^0, \Phi_f^1, \Phi_f^2\}$ be the nodal values of the level set function $\Phi_f$. In this case, we compute the right-hand side for side $a$ to get the equivalent polynomial $H_{eq}^a$, and similarly for side $b$:
\begin{align}
    \label{eq:regular-cut-cell-ab}
    b^a_j &= \int_{K_r} H(\Phi_f) \varphi_j \, dx = \int_{K_a} \varphi_j \, dx , \\
    b^b_j &= \int_{K_r} (1-H(\Phi_f)) \varphi_j \, dx = \int_{K_b} \varphi_j \, dx .
\end{align}

Although this approach has been used in the context of CutFEM for problems with a single interface~\cite{keesSimpleAccurateEfficient2022}, we are not aware of any work extending it to multiple interfaces. One of the challenges in unfitted problems with two or more intersecting interfaces is integrating over the cells that are cut by more than one interface.
\subsection{Triply-cut cell}
\begin{figure}[!hb]
    \centering
    \includegraphics[width=0.95\linewidth]{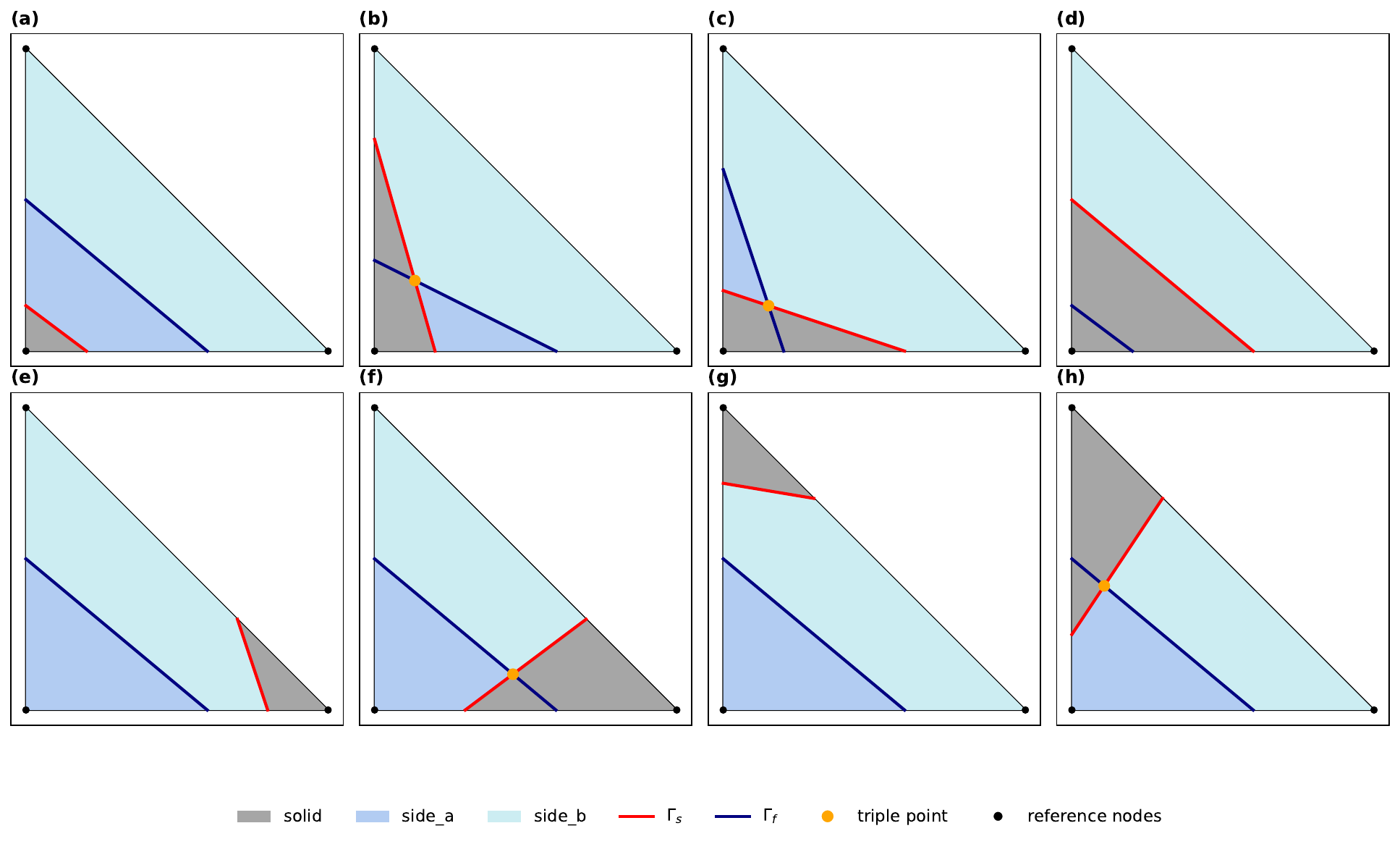}
    \caption{The distinct topologies of a triply-cut cell: (a)--(d) the four same-node
    orderings of $K_s$ and $K_a$ (nested, the two crossing/triple-point cases, and
    reversed); (e)--(h) the four configurations where $\Gamma_s$ isolates a different
    node than $\Gamma_f$, each with and without a triple point.}
    \label{fig:triply-cut-cell-configs}
\end{figure}

Consider a cell $K_r \in \mathcal{T}_h$ that is cut by both interfaces $\Gamma_s$ and $\Gamma_f$, so that it carries two independent sets of nodal level set values, $\{\Phi_s^0,\Phi_s^1,\Phi_s^2\}$ and $\{\Phi_f^0,\Phi_f^1,\Phi_f^2\}$. As shown in Fig.~\ref{fig:triangulation-schematic}, such a cell is split into three regions: the solid and the two fluid phases $a$ and $b$. As in the typical FEM machinery, the integrals are computed on a reference element. We map the corner of the element that is isolated by $\Gamma_f$ to node $0$ of the reference triangle, and, for the description below, call the phase of this corner $a$; when it is phase $b$, the roles of $a$ and $b$ are exchanged. This choice makes it easier to lay out all the possible configurations of the two interface cuts. Fig.~\ref{fig:triply-cut-cell-configs} shows the eight possible configurations of the cuts, categorized by which node is isolated by $\Gamma_s$ and whether the two cuts cross to produce a triple point inside the element or not. The four configurations in panels (a)--(d) have $\Gamma_s$ isolating the same node as $\Gamma_f$, while the four configurations in panels (e)--(h) have $\Gamma_s$ isolating a different node. Using this categorization, we can compute the right-hand side vector $\textbf{b}$ for each of the three regions in the triply-cut cell.

Unlike the regular cut cell, computing $\textbf{b}$ here is not a single moment-fitting solve. Each phase region has its own equivalent polynomial, computed from that phase's right-hand side $\textbf{b}$. By construction
\begin{align}
    b_j^{solid} + b_j^{a} + b_j^{b} \equiv \int_{K_r} \varphi_j \, dx ,
\end{align}
and the equation to solve for the equivalent polynomial of each phase is
\begin{align}
    \mathbb{M} \hat{\textbf{H}}^{phase} = \textbf{b}^{phase}, \qquad \text{phase} \in \{\text{solid}, a, b\}.
\end{align}
The solid right-hand side is exactly the regular single-interface formula applied to $\Gamma_s$,
\begin{align}
    b_j^{solid} = \int_{K_r} (1-H(\Phi_s)) \varphi_j \, dx = \int_{K_s} \varphi_j \, dx.
\end{align}
The right-hand sides of phases $a$ and $b$, however, depend on the cut configuration explained above. For the case where $\Gamma_s$ cuts the same corner as $\Gamma_f$, there are three cases, distinguished by comparing the two pairs of cut points and shown in panels (a)--(d) of Fig.~\ref{fig:triply-cut-cell-configs}:
\begin{itemize}
    \item \textbf{Nested} (Fig.~\ref{fig:triply-cut-cell-configs}a), $x_s \le x_f$ and $y_s \le y_f$ ($K_s \subseteq K_a$): the solid corner sits entirely inside the corner of phase $a$, so phase $a$ is simply $K_a$ with the solid removed, and phase $b$ is untouched,
    \begin{align}
        b_j^{a} = \int_{K_a} \varphi_j \, dx - \int_{K_s} \varphi_j \, dx, \qquad
        b_j^{b} = \int_{K_r} \varphi_j \, dx - b_j^{solid} - b_j^{a},
    \end{align}
    where $\int_{K_a} \varphi_j \, dx$ is the ordinary fluid--fluid moment of Eq.~\eqref{eq:regular-cut-cell-ab} computed on $K_a$ alone.
    \item \textbf{Reversed} (Fig.~\ref{fig:triply-cut-cell-configs}d), $x_s \ge x_f$ and $y_s \ge y_f$ ($K_a \subseteq K_s$): $\Gamma_f$ lies entirely inside the solid. This is really an ordinary cell cut only by $\Gamma_s$. Hence, the fluid right-hand side (phase $b$) is simply
    \begin{align}
        b_j^{b} = \int_{K_r} \varphi_j \, dx - b_j^{solid}.
    \end{align}
    \item \textbf{Crossing} (Fig.~\ref{fig:triply-cut-cell-configs}b,~\ref{fig:triply-cut-cell-configs}c), $x_s < x_f, y_s > y_f$ (or its mirror $x_s > x_f, y_s < y_f$): the two cuts cross at a genuine triple point $P$, so part of $K_a$ lies inside the solid $K_s$. The overlap can be located exactly using the triple point $P$. Writing the two cut lines as $x/x_s+y/y_s=1$ and $x/x_f+y/y_f=1$, they intersect at
    \begin{align}
        P = \left( \frac{(y_s-y_f)\,x_s x_f}{x_f y_s - x_s y_f}, \ \frac{(x_f-x_s)\,y_s y_f}{x_f y_s - x_s y_f} \right),
    \end{align}
    and the true region of phase $a$ is the single triangle $(x_s,0),(x_f,0),P$ when $x_s<x_f,y_s>y_f$, or its mirror $(0,y_f),(0,y_s),P$ when $x_s>x_f,y_s<y_f$. The right-hand side of phase $a$ is then computed directly on that triangle,
    \begin{align}
        b_j^{a} = \int_{\triangle(x_s,0),(x_f,0),P} \varphi_j \, dx \qquad \text{(or its mirror)},
    \end{align}
    and phase $b$ is again obtained by complement, $b_j^{b} = \int_{K_r}\varphi_j\,dx - b_j^{solid} - b_j^{a}$.
\end{itemize}

For the case where $\Gamma_s$ cuts a different corner than $\Gamma_f$, shown in panels (e)--(h) of Fig.~\ref{fig:triply-cut-cell-configs}, the formula of Eq.~\eqref{eq:cut-formula} also gives the cut point of $\Gamma_s$ located on the hypotenuse. For an edge from node $i$ to node $j$ with $\Phi_s^i,\Phi_s^j$ of opposite sign, the zero of the affine interpolant of $\Phi_s$ along it, as a fraction of the way from node $i$ to node $j$, is
\begin{align}
    \label{eq:hypotenuse-cut}
    z = \frac12 - \frac12\frac{\Phi_s^j+\Phi_s^i}{\Phi_s^j-\Phi_s^i}.
\end{align}
If $\Gamma_s$ isolates node $1$, we take $(i,j)=(1,2)$, and the cut points are $(x_s,0)$ and $(1-z,z)$. If instead $\Gamma_s$ isolates node $2$, we take $(i,j)=(2,1)$, so that $z$ is measured from node $2$, and the cut points are $(z,1-z)$ and $(0,y_s)$; the two cases are mirror images under the exchange of $x$ and $y$. With this cut point calculation, we can compute the right-hand side both when the two cuts produce a triple point and when they do not, covering all four cases.

\begin{itemize}
    \item \textbf{No triple point} (Fig.~\ref{fig:triply-cut-cell-configs}e,~\ref{fig:triply-cut-cell-configs}g): $K_s$ and $K_a$ are disjoint corner triangles sitting at different vertices, so the right-hand side for each phase is just the regular single-interface formula applied to that phase's triangle,
    \begin{align}
        b_j^{a} = \int_{K_a} \varphi_j \, dx, \qquad b_j^{b} = \int_{K_r}\varphi_j\,dx - b_j^{solid} - b_j^{a}.
    \end{align}
    \item \textbf{Triple point} (Fig.~\ref{fig:triply-cut-cell-configs}f,~\ref{fig:triply-cut-cell-configs}h): computing the vector $\textbf{b}$ requires the triple point $P$ in order to locate the regions of each phase.
        \begin{itemize}
            \item If $\Gamma_s$ isolates node $1$, the two cut lines are $\Gamma_f$, $x/x_f+y/y_f=1$, and $\Gamma_s$, $z x - (1 - z - x_s) y = z x_s$, with the intersection point
                \begin{align}
                    P = \left( \frac{x_f\big[z x_s + (1-z-x_s)\,y_f\big]}{z x_f + (1-z-x_s)\,y_f}, \ \frac{z\,y_f\,(x_f-x_s)}{z x_f + (1-z-x_s)\,y_f} \right).
                \end{align}
                The part of $K_a$ that overlaps $K_s$ is itself the triangle with vertices $(x_s,0)$, $(x_f,0)$ and $P$. Removing it from $K_a$ leaves a \emph{quadrilateral} with vertices $(0,0)$, $(x_s,0)$, $P$ and $(0,y_f)$. The right-hand side for this quadrilateral is computed by splitting it into two triangles,
                \begin{align}
                    b_j^{a} = \int_{\triangle(0,0),(x_s,0),P} \varphi_j \, dx + \int_{\triangle(0,0),P,(0,y_f)} \varphi_j \, dx ,
                \end{align}
                and phase $b$ is again obtained by complement, $b_j^{b} = \int_{K_r}\varphi_j\,dx - b_j^{solid} - b_j^{a}$.
            \item If $\Gamma_s$ isolates node $2$, the two cut lines are $\Gamma_f$, $x/x_f+y/y_f=1$, and $\Gamma_s$, $z y - (1 - z - y_s) x = z y_s$, with the intersection point
                \begin{align}
                    P = \left( \frac{z\,x_f\,(y_f-y_s)}{z y_f + (1-z-y_s)\,x_f}, \ \frac{y_f\big[z y_s + (1-z-y_s)\,x_f\big]}{z y_f + (1-z-y_s)\,x_f} \right).
                \end{align}
                The right-hand side is computed in the same way as above, with the quadrilateral with vertices $(0,0)$, $(0,y_s)$, $P$ and $(x_f,0)$ split into two triangles.
        \end{itemize}
\end{itemize}

\section{Results}
\label{sec:results}

We test the method on a three-phase problem with a straight fluid--fluid interface and a circular solid. We report the convergence of linear ($p=1$) and quadratic ($p=2$) elements and compare the method with the standard finite element method on a mesh fitted to both interfaces.

\begin{figure}[!ht]
    \centering
    \includegraphics[width=0.5\linewidth]{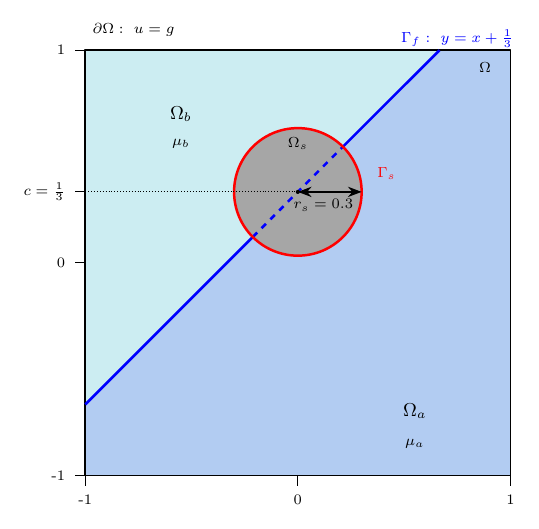}
    \caption{The geometry of the test problem.}
    \label{fig:test-geometry}
\end{figure}

The domain is $\Omega = [-1,1]^2$. The solid is a disk of radius $r=0.3$ centered at $(0,c)$, and $\Gamma_f$ is the line $y = mx + c$ with $m=1$ and $c=\frac13$, so that $\Gamma_f$ passes through the center of the disk. The exact geometry of the test problem is illustrated in Fig.~\ref{fig:test-geometry}. The level set functions are
\begin{align}
    \Phi_s(x,y) &= \sqrt{x^2 + (y-c)^2} - r, \\
    \Phi_f(x,y) &= \frac{m x + c - y}{\sqrt{m^2 + 1}},
\end{align}
so that $\Omega_a$ lies below $\Gamma_f$ and $\Omega_b$ above it.

The penalty parameters used for all tests are as follows: 1) Nitsche parameter, $\gamma_N = 10\,p^2$, since the constant of the trace inequality grows like $p^2$~\cite{warburtonConstantsHpfiniteElement2003}, 2) Ghost penalty parameter, $\gamma_g = 0.15$, and 3) Interior penalty parameter, $\gamma_I=100$ (only for $p=2$). The unfitted method uses Gauss quadrature of order 6 for $p=1$ and order 8 for $p=2$ on the elements and on the element edges, so that the equivalent polynomials and the Nitsche and edge terms are integrated accurately; the mesh fitted to both interfaces uses order 6. The errors are computed over the fluid domain.

\textit{Test 1:} The manufactured solution is
\begin{align}
    \label{eq:test1-mms}
    u = \begin{cases}
        \dfrac{g}{\mu_a} & \text{in } \Omega_a, \\[2ex]
        \dfrac{g}{\mu_b} & \text{in } \Omega_b,
    \end{cases} \qquad
    g = \sin(\pi \xi) \cos(\pi \eta), \qquad f = 2\pi^2 g .
\end{align}
Here $\xi$ and $\eta$ are the orthonormal coordinates attached to $\Gamma_f$; $\xi$ being the normal coordinate and $\eta$ being the tangential coordinate with the origin at $(0,c)$.
\begin{align}
    \xi &= \frac{y - m x - c}{\sqrt{m^2 + 1}}, \quad \eta = \frac{x + m (y - c)}{\sqrt{m^2 + 1}}
\end{align}
\begin{figure}[!ht]
    \centering
    \begin{subfigure}{0.45\linewidth}
        \includegraphics[width=\linewidth]{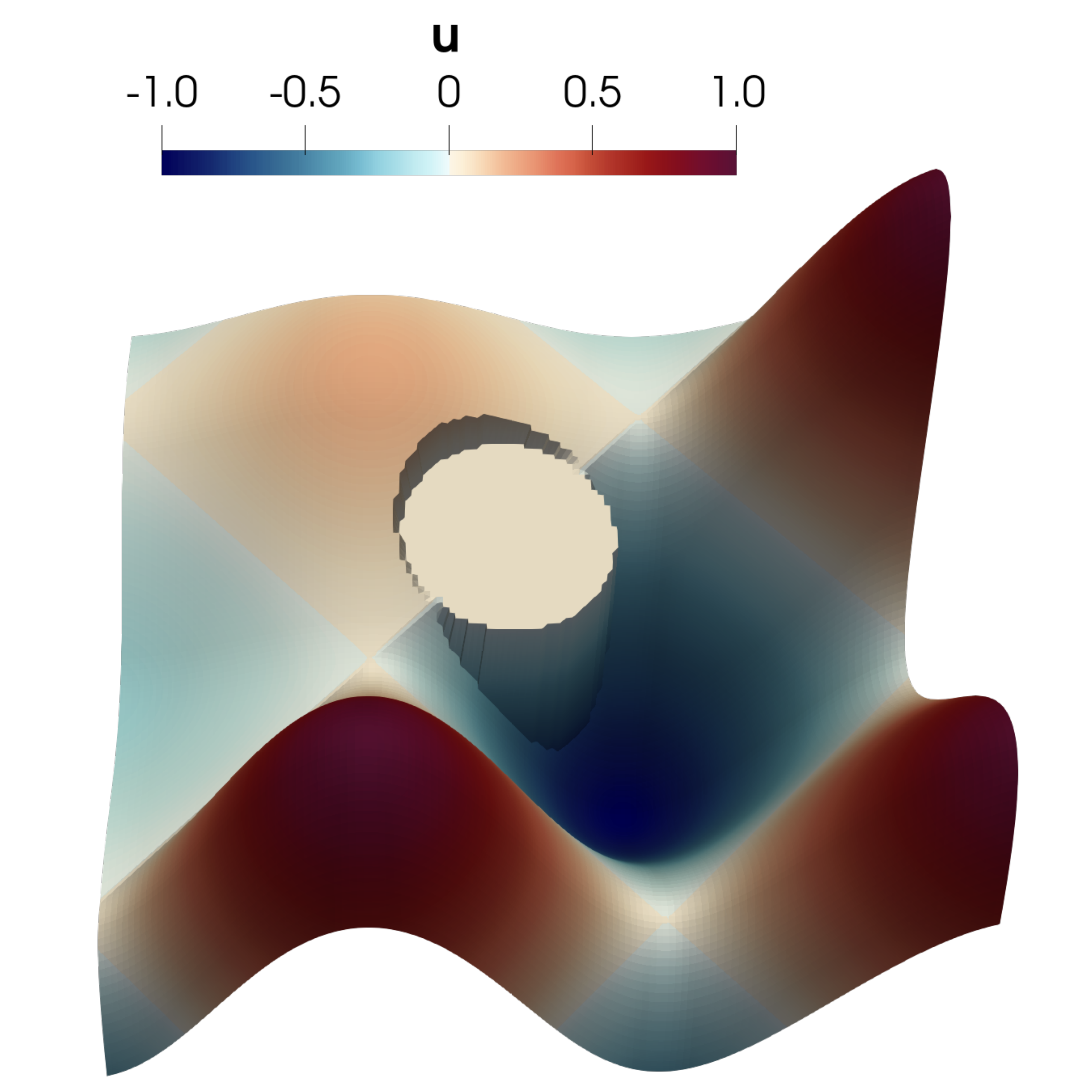}
        \caption{Test 1}
        \label{fig:test1-solution}
    \end{subfigure}
    \hfill
    \begin{subfigure}{0.45\linewidth}
        \includegraphics[width=\linewidth]{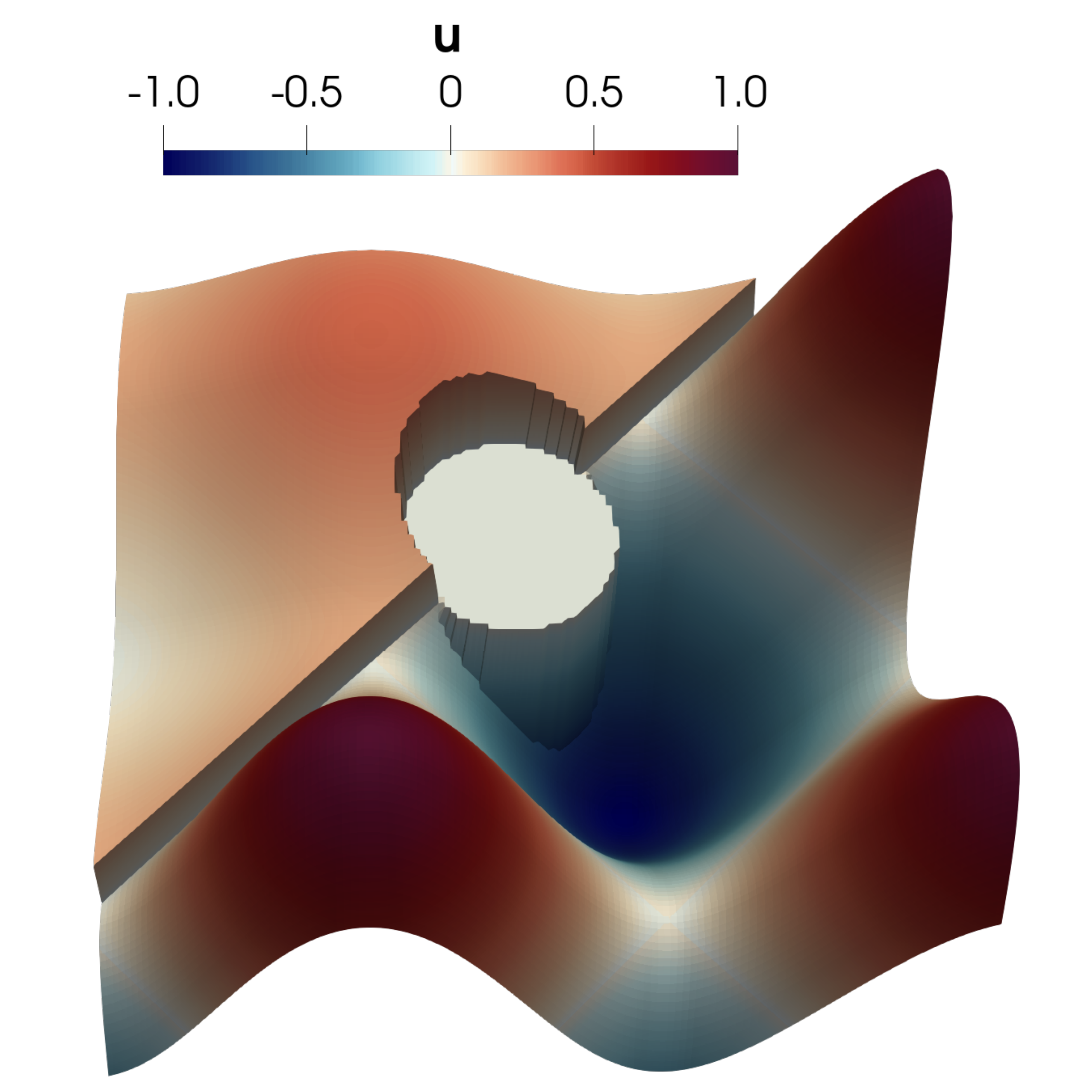}
        \caption{Test 2 (with $C = 0.2$)}
        \label{fig:test2-solution}
    \end{subfigure}
    \caption{Solution plots for $\mu_b/\mu_a = 5$.}
    \label{fig:solution}
\end{figure}

By performing this coordinate transformation, the manufactured solution is exactly zero on $\Gamma_f$ (regardless of $m$ and $c$) and changes in the normal direction $\xi$ away from the interface. Note that $\xi$ is exactly the signed distance from the fluid--fluid interface. A crucial point to note is that the orthonormal coordinates are just a rigid translation and rotation and therefore the Laplacian is invariant. This ensures that the Poisson equation retains its form under the coordinate transformation. The jump conditions across the fluid--fluid interface are therefore preserved in the new coordinate system. Fig.~\ref{fig:test1-solution} shows the solution plot.

First, we compare the presented unfitted method with the mesh fitted to both interfaces. The fitted mesh is an unstructured mesh, while both structured and unstructured meshes are used for the unfitted method for comparison. Table~\ref{tab:fitted} compares the $L^2$ errors at about the same number of nodes. All three discretizations converge at the optimal order $p+1$. The unfitted method thus attains the accuracy of a fitted mesh with a single background mesh for any position of the interfaces, whereas the fitted mesh has to be regenerated whenever an interface moves and must contain the triple points as vertices.
\begin{table}[!ht]
    \centering
    \footnotesize
    \renewcommand{\arraystretch}{1.3}
    \setlength{\tabcolsep}{3pt}
    \caption{Test~1, $\mu_b/\mu_a=5$: $L^2$ error and rate on the fitted mesh and for the unfitted method on unstructured and structured meshes, against the number of nodes $N$.}
        \begin{tabular}{c||cc|cc|cc}
        & \multicolumn{2}{c|}{fitted mesh} & \multicolumn{2}{c|}{unfitted, unstructured} & \multicolumn{2}{c}{unfitted, structured} \\
        $N$ & $\|u^h-u\|_{L^2}$ & rate & $\|u^h-u\|_{L^2}$ & rate & $\|u^h-u\|_{L^2}$ & rate \\\hline\hline
        \multicolumn{7}{c}{\textbf{$p=1$}} \\\hline
        $\sim100$ & $6.415\times10^{-2}$ & --- & $5.937\times10^{-2}$ & --- & $5.654\times10^{-2}$ & --- \\
        $\sim400$ & $1.638\times10^{-2}$ & 1.97 & $1.370\times10^{-2}$ & 2.12 & $1.377\times10^{-2}$ & 2.04 \\
        $\sim1\,600$ & $3.787\times10^{-3}$ & 2.11 & $3.777\times10^{-3}$ & 1.86 & $3.404\times10^{-3}$ & 2.02 \\
        $\sim6\,400$ & $9.753\times10^{-4}$ & 1.96 & $9.582\times10^{-4}$ & 1.98 & $8.502\times10^{-4}$ & 2.00 \\
        $\sim25\,600$ & $2.476\times10^{-4}$ & 1.98 & $2.412\times10^{-4}$ & 1.99 & $2.125\times10^{-4}$ & 2.00 \\
        $\sim102\,400$ & $6.202\times10^{-5}$ & 2.00 & $6.100\times10^{-5}$ & 1.98 & $5.321\times10^{-5}$ & 2.00 \\\hline\hline
        \multicolumn{7}{c}{\textbf{$p=2$}} \\\hline
        $\sim400$ & $3.212\times10^{-3}$ & --- & $3.186\times10^{-3}$ & --- & $3.152\times10^{-3}$ & --- \\
        $\sim1\,600$ & $4.086\times10^{-4}$ & 2.97 & $3.990\times10^{-4}$ & 3.00 & $3.993\times10^{-4}$ & 2.98 \\
        $\sim6\,400$ & $5.354\times10^{-5}$ & 2.93 & $5.236\times10^{-5}$ & 2.93 & $4.930\times10^{-5}$ & 3.02 \\
        $\sim25\,600$ & $6.600\times10^{-6}$ & 3.02 & $6.816\times10^{-6}$ & 2.94 & $6.194\times10^{-6}$ & 2.99 \\
        $\sim102\,400$ & $8.558\times10^{-7}$ & 2.95 & $8.464\times10^{-7}$ & 3.01 & $7.717\times10^{-7}$ & 3.00 \\
        $\sim409\,600$ & $1.061\times10^{-7}$ & 3.01 & $1.071\times10^{-7}$ & 2.98 & $9.657\times10^{-8}$ & 3.00
    \end{tabular}
    \label{tab:fitted}
\end{table}

Table~\ref{tab:test1-rate-1000} shows Test~1 with $\mu_b/\mu_a = 1000$. For $p=1$ the $L^2$ error converges at second order while for $p=2$ the $L^2$ error converges at third order. For both polynomial orders, the $L^\infty$ error fluctuates around order 2 for $p=1$ and 3 for $p=2$. The fluctuations for high coefficient contrast are due to the large discontinuity of the IFEM basis on the cut edges. Without using the interior penalty, the $L^\infty$ errors degrade severely.

\begin{table}[htbp]
    \centering
    \footnotesize
    \renewcommand{\arraystretch}{1.3}
    \setlength{\tabcolsep}{3.5pt}
    \caption{Test~1 with the high contrast $\mu_a=1,\mu_b=1000$ on the structured meshes, $p=1$ and $p=2$.}
    \begin{tabular}{c|cc|cc||cc|cc}
        & \multicolumn{4}{c||}{$p=1$} & \multicolumn{4}{c}{$p=2$} \\
        $h$ & $\|u^h-u\|_{L^2}$ & rate & $\|u^h-u\|_{L^\infty}$ & rate & $\|u^h-u\|_{L^2}$ & rate & $\|u^h-u\|_{L^\infty}$ & rate \\\hline
        $1/4$ & $9.236\times10^{-2}$ & --- & $1.491\times10^{-1}$ & --- & $7.429\times10^{-3}$ & --- & $2.067\times10^{-2}$ & --- \\
        $1/8$ & $2.187\times10^{-2}$ & 2.08 & $3.616\times10^{-2}$ & 2.04 & $8.826\times10^{-4}$ & 3.07 & $3.868\times10^{-3}$ & 2.42 \\
        $1/16$ & $5.355\times10^{-3}$ & 2.03 & $8.878\times10^{-3}$ & 2.03 & $9.797\times10^{-5}$ & 3.17 & $3.188\times10^{-4}$ & 3.60 \\
        $1/32$ & $1.323\times10^{-3}$ & 2.02 & $3.392\times10^{-3}$ & 1.39 & $1.235\times10^{-5}$ & 2.99 & $6.391\times10^{-5}$ & 2.32 \\
        $1/64$ & $3.295\times10^{-4}$ & 2.01 & $5.439\times10^{-4}$ & 2.64 & $1.501\times10^{-6}$ & 3.04 & $5.175\times10^{-6}$ & 3.63 \\
        $1/128$ & $8.237\times10^{-5}$ & 2.00 & $1.356\times10^{-4}$ & 2.00 & $1.882\times10^{-7}$ & 3.00 & $1.002\times10^{-6}$ & 2.37
    \end{tabular}
    \label{tab:test1-rate-1000}
\end{table}

\noindent
\textit{Test 2:} This test adds a constant solution jump across $\Gamma_f$ to Test~1:
\begin{align}
    \label{eq:test2-mms}
    u =
    \begin{cases}
        \dfrac{g}{\mu_a} & \text{in } \Omega_a, \\[2ex]
        \dfrac{g}{\mu_b} + C & \text{in } \Omega_b,
    \end{cases}
\end{align}
where $C$ is a constant, so that $J_s = C$; the flux and the source are those of Test~1. The solution is shown in Fig.~\ref{fig:test2-solution}.

\vspace{5mm}
\noindent
\textit{Test 3:} This test adds a constant flux jump across $\Gamma_f$ to Test~1:
\begin{align}
    \label{eq:test3-mms}
    u =
    \begin{cases}
        \dfrac{g}{\mu_a} & \text{in } \Omega_a, \\[2ex]
        \dfrac{g + C \xi}{\mu_b} & \text{in } \Omega_b,
    \end{cases}
\end{align}
where $C$ is a constant, so that $J_f=C$.

\begin{table}[!ht]
    \centering
    \footnotesize
    \renewcommand{\arraystretch}{1.25}
    \caption{Tests 2 and 3 on the structured meshes, $\mu_a=1,\mu_b=5$, $p=1$ and $p=2$. At each mesh size and for each norm, the two tests are printed with the number of significant digits at which they first differ.}
    \begin{tabular}{c|c|c|c|c}
        $h$ & $\|u^h-u\|_{L^2(\Omega_f)}$ & rate & $\|u^h-u\|_{L^\infty(\Omega_f)}$ & rate \\\hline\hline
        \multicolumn{5}{c}{\textbf{Test 2: solution jump $C=0.2$, $p=1$}} \\\hline
        $1/4$ & $9.107175999454194\times10^{-2}$ & --- & $1.452\times10^{-1}$ & --- \\
        $1/8$ & $2.16608928738921\times10^{-2}$ & 2.07 & $3.5536516762443\times10^{-2}$ & 2.03 \\
        $1/16$ & $5.3593741479202\times10^{-3}$ & 2.01 & $8.8283635680410\times10^{-3}$ & 2.01 \\
        $1/32$ & $1.3272844391412\times10^{-3}$ & 2.01 & $2.182\times10^{-3}$ & 2.02 \\
        $1/64$ & $3.32327857219\times10^{-4}$ & 2.00 & $5.44020599276\times10^{-4}$ & 2.00 \\
        $1/128$ & $8.31077469078\times10^{-5}$ & 2.00 & $1.35642570826\times10^{-4}$ & 2.00 \\\hline
        \multicolumn{5}{c}{\textbf{Test 3: flux jump $C=0.2$, $p=1$}} \\\hline
        $1/4$ & $9.107175999454192\times10^{-2}$ & --- & $1.452\times10^{-1}$ & --- \\
        $1/8$ & $2.16608928738922\times10^{-2}$ & 2.07 & $3.5536516762444\times10^{-2}$ & 2.03 \\
        $1/16$ & $5.3593741479204\times10^{-3}$ & 2.01 & $8.8283635680413\times10^{-3}$ & 2.01 \\
        $1/32$ & $1.3272844391414\times10^{-3}$ & 2.01 & $2.182\times10^{-3}$ & 2.02 \\
        $1/64$ & $3.32327857220\times10^{-4}$ & 2.00 & $5.44020599277\times10^{-4}$ & 2.00 \\
        $1/128$ & $8.31077469083\times10^{-5}$ & 2.00 & $1.35642570827\times10^{-4}$ & 2.00 \\\hline
        \multicolumn{5}{c}{\textbf{Test 2: solution jump $C=0.2$, $p=2$}} \\\hline
        $1/4$ & $6.48027009321061\times10^{-3}$ & --- & $1.27272815338\times10^{-2}$ & --- \\
        $1/8$ & $7.67619227128\times10^{-4}$ & 3.08 & $1.37763887632\times10^{-3}$ & 3.21 \\
        $1/16$ & $9.71466899882\times10^{-5}$ & 2.98 & $2.2100360011\times10^{-4}$ & 2.64 \\
        $1/32$ & $1.20489984109\times10^{-5}$ & 3.01 & $2.283690871\times10^{-5}$ & 3.27 \\
        $1/64$ & $1.5104474861\times10^{-6}$ & 3.00 & $3.5053441\times10^{-6}$ & 2.70 \\
        $1/128$ & $1.884307869\times10^{-7}$ & 3.00 & $3.591840\times10^{-7}$ & 3.29 \\\hline
        \multicolumn{5}{c}{\textbf{Test 3: flux jump $C=0.2$, $p=2$}} \\\hline
        $1/4$ & $6.48027009321062\times10^{-3}$ & --- & $1.27272815339\times10^{-2}$ & --- \\
        $1/8$ & $7.67619227127\times10^{-4}$ & 3.08 & $1.37763887630\times10^{-3}$ & 3.21 \\
        $1/16$ & $9.71466899879\times10^{-5}$ & 2.98 & $2.2100360014\times10^{-4}$ & 2.64 \\
        $1/32$ & $1.20489984108\times10^{-5}$ & 3.01 & $2.283690873\times10^{-5}$ & 3.27 \\
        $1/64$ & $1.5104474858\times10^{-6}$ & 3.00 & $3.5053442\times10^{-6}$ & 2.70 \\
        $1/128$ & $1.884307866\times10^{-7}$ & 3.00 & $3.591844\times10^{-7}$ & 3.29
    \end{tabular}
    \label{tab:tests123}
\end{table}

Table~\ref{tab:tests123} shows optimal convergence, of order $p+1$, in both norms and for both degrees. For $p=2$ the $L^\infty$ rate alternates between $2.6$ and $3.3$, with an average of $3.0$. The errors of both tests equal those of Test~1 on the same meshes up to round-off because $j_h$ and $q_h$ reproduce exactly the constant solution jump and flux jump. For the same reason the errors do not depend on the size of the jumps.

\vspace{5mm}
\noindent
\textit{Test 4:} This test replaces the Dirichlet condition on $\Gamma_s$ by a Robin condition, which models convective heat transfer (Newton's law of cooling):
\begin{align}
    \label{eq:robin}
    -\mu \frac{\partial u}{\partial n_s} = h_c (u - u_\infty) \qquad \text{on } \Gamma_s,
\end{align}
where $h_c > 0$ is the heat transfer coefficient and $u_\infty$ the ambient temperature. Since $\Gamma_s$ borders both fluids, the coefficients $h_c^a$ and $h_c^b$ are used on the parts of $\Gamma_s$ bordering $\Omega_a$ and $\Omega_b$. The Robin condition replaces the boundary flux term of the weak form,
\begin{align*}
    -\int_{\Gamma_{s}} \mu \frac{\partial u^h}{\partial n_s} v^h \, ds = \int_{\Gamma_s} h_c\,(u^h - u_\infty)\, v^h\, ds ,
\end{align*}
so that the Nitsche terms of $a$ and $l$ on $\Gamma_s$ are replaced by
\begin{align}
    \label{eq:robin-terms}
    a_{\Gamma_s}(u^h,v^h) = \int_{\Gamma_{s}} h_c\, u^h\, v^h\, ds , \qquad
    l_{\Gamma_s}(v^h) = \int_{\Gamma_{s}} h_c\, u_\infty\, v^h\, ds .
\end{align}
The exact solution is that of Test~1 with $\mu_a = 1$ and $\mu_b = 5$, the coefficients are $h_c^a = 600$ and $h_c^b = 15$, and $u_\infty = u + \mu\,(\partial u/\partial n_s)/h_c$ is manufactured from the exact solution.

\begin{table}[!ht]
    \centering
    \footnotesize
    \renewcommand{\arraystretch}{1.3}
    \setlength{\tabcolsep}{3.5pt}
    \caption{Test 4 with a Robin condition on the fluid--solid interface $\Gamma_s$ ($\mu_a=1,\mu_b=5$, $h_c^a=600$, $h_c^b=15$).}
    \begin{tabular}{c|cc|cc||cc|cc}
        & \multicolumn{4}{c||}{$p=1$} & \multicolumn{4}{c}{$p=2$} \\
        $h$ & $\|u^h-u\|_{L^2}$ & rate & $\|u^h-u\|_{L^\infty}$ & rate & $\|u^h-u\|_{L^2}$ & rate & $\|u^h-u\|_{L^\infty}$ & rate \\\hline
        $1/4$ & $9.177\times10^{-2}$ & --- & $1.498\times10^{-1}$ & --- & $6.596\times10^{-3}$ & --- & $1.281\times10^{-2}$ & --- \\
        $1/8$ & $2.202\times10^{-2}$ & 2.06 & $3.579\times10^{-2}$ & 2.07 & $7.755\times10^{-4}$ & 3.09 & $1.359\times10^{-3}$ & 3.24 \\
        $1/16$ & $5.415\times10^{-3}$ & 2.02 & $8.855\times10^{-3}$ & 2.01 & $1.066\times10^{-4}$ & 2.86 & $2.331\times10^{-4}$ & 2.54 \\
        $1/32$ & $1.336\times10^{-3}$ & 2.02 & $2.184\times10^{-3}$ & 2.02 & $1.535\times10^{-5}$ & 2.80 & $4.884\times10^{-5}$ & 2.26 \\
        $1/64$ & $3.342\times10^{-4}$ & 2.00 & $5.449\times10^{-4}$ & 2.00 & $2.888\times10^{-6}$ & 2.41 & $1.282\times10^{-5}$ & 1.93 \\
        $1/128$ & $8.345\times10^{-5}$ & 2.00 & $1.358\times10^{-4}$ & 2.01 & $6.267\times10^{-7}$ & 2.20 & $3.683\times10^{-6}$ & 1.80
    \end{tabular}
    \label{tab:test4-rate}
\end{table}

Table~\ref{tab:test4-rate} shows second-order convergence for $p=1$. For $p=2$ the rate decreases toward two. The Robin condition depends on the normal of $\Gamma_s$, and the normal of the chord $\Gamma_{s,h}$ differs from that of the circle by $O(h)$. To reach optimal convergence for $p=2$, a higher order level-set representation is required, which we leave for future work.

\section{Conclusion}
\label{sec:conclusion}

We proposed a hybrid CutFEM--IFEM method for the three-phase Poisson problem, which combines a Nitsche-based CutFEM treatment of the embedded fluid--solid interface with an IFEM treatment of the immersed fluid--fluid interface, for linear and quadratic elements. Cells cut by both interfaces are integrated with composite equivalent polynomials for the Heaviside and Dirac delta functions. The IFEM basis functions are enriched with the solution and flux jump lift functions to represent non-homogeneous jumps. The enrichment functions reproduce the jumps exactly so that the error does not depend on the size of the jumps. The numerical experiments show optimal convergence of order $p+1$ for both degrees, with errors comparable to those of a mesh fitted to both interfaces. Furthermore, the optimal convergence is observed for the coefficient contrast up to $1000$. With a Robin condition on a curved boundary, the quadratic rate is limited to two by the linear representation of the boundary. A quadratic basis that is bounded independently of the coefficient contrast, and a higher-order representation of the embedded boundary, are the subject of future work.
\section*{Code availability}

The method is implemented in the Proteus toolkit. The code used for all results in this paper is the tag \path{hybrid_cutfem_ifem-v1} of the repository \url{https://github.com/cekees/proteus} (branch \path{3phase-hybrid-cutfem-ifem}, based on the Proteus 1.9.x line). The drivers are \path{test/ci/ladr_3p_2d_p.py} and \path{test/ci/ladr_3p_2d_c0pk_n.py} for the unfitted method and \path{test/ci/ladr_fitted_2d_p.py} and \path{test/ci/ladr_fitted_2d_n.py} for the mesh fitted to both interfaces. At this tag the runs also require the environment variables \path{ADR_DISABLE_TC=1}, \path{ADR_TC_EDGE_FLUID=1} and \path{ADR_TC_EDGE_BOUNDARY_ONLY=1}, and the penalty parameters and quadrature orders given in Section~\ref{sec:results}.

\bibliographystyle{elsarticle-num}
\bibliography{main}

\end{document}